\documentclass[titlepage,11pt,a4paper]{article}
\usepackage{amstext}
\usepackage{amssymb}
\usepackage{amsfonts}
\usepackage{amssymb}
\usepackage{amsmath}
\usepackage{graphicx}
\usepackage[arrow,matrix,curve]{xy}
\usepackage[latin1]{inputenc}
\usepackage{theorem}
\newtheorem{satz}{Theorem}[section]
\newtheorem{defi}[satz]{Definition}
\newtheorem{bem}[satz]{Remark}
\newtheorem{kor}[satz]{Corollary}
\newtheorem{lem}[satz]{Lemma}

\newtheorem{pro}[satz]{Proposition}

\newcommand{\modu}{~\rm{mod}~}

\newcommand{\qed}{\begin{flushright}$\square$\end{flushright}}
\newcommand{\filt}[6]{
\[
\begin{xy}
\xymatrix@R20pt@C20pt{
&\mathbb{C}^3&\\\langle #1,#2 \rangle\ar[ru]&\langle #3,#4\rangle\ar[u]&\langle #5,#6\rangle\ar[lu]\\\langle #1\rangle\ar[u]&\langle #3\rangle\ar[u]&\langle #5\rangle\ar[u]}
\end{xy}
\]
} 
\begin{document}
\begin{center} 
{\bf\LARGE Torus fixed points of moduli spaces of
stable bundles of rank three
\\}
\end{center}
\vspace{0.0cm}
\begin{center}
Thorsten Weist\\Fachbereich C - Mathematik\\
Bergische Universität Wuppertal\\
D - 42097 Wuppertal, Germany\\
e-mail: weist@math.uni-wuppertal.de
\end{center}
\vspace{0.0cm}
\section{Introduction}
The aim of this paper is to extend the main result of \cite{kly1}, which determines
the Euler characteristic of moduli spaces of stable bundles of rank two on
the projective plane, to the case of bundles of rank three. Therefore, we first
discuss some results of \cite{kly1}, \cite{kly2} and \cite{kly3} and some general results
about vector bundles and their moduli spaces. Klyachko proved that toric
bundles of rank $n$ correspond to filtrations of an $n$-dimensional vector space.
These filtrations can also be understood as representations of the subspace
quiver. Since the stability condition can be transferred, the moduli spaces
of stable representations can be identified with some fixed point components
of the moduli space of bundles. Toric bundles of rank three correspond to
filtrations of a three dimensional vector space.\\
The length of the arms of the subspace quiver can be used to determine
the Chern classes of the corresponding torus fixed points. Investigating the
moduli spaces of (semi-)stable representations it turns out that the second
Chern class varies depending on how many two-dimensional subspaces the
one-dimensional ones contain. For the discriminant of bundles of rank three
on the projective plane we either have 0 or 4 $\modu$ 6. In the first case we
first analyse the polystable points which are representations which can be
decomposed into representations of the same slope. It turns out that the
existence and number of those points just depend on the length of the arms of
the subspace quiver. We also investigate those stable points of the moduli
space having different second Chern class. In the second case there exist
only stable points, i.e. there are no polystable points, so we just have to
consider the latter cases. In both cases it turns out that the moduli spaces
of semistable representations are projective lines and that there exist only
finitely many of the described points so that it is easy to calculate the Euler
characteristic of the fixed point components.\\
The stability condition reduces to a system of linear inequalities which can
be solved using standard methods from e.g. \cite{pad}. After stating how to
determine the solutions of such a system we will see that in our case they                                                                                                                                                                                                        
correspond to solutions of some quadratic equations so that the generating
functions of the Euler characteristic of the moduli spaces of stable bundles
of rank three on the projective plane can be determined.\\\\
{\bf Acknowledgment:} I would like to thank Markus Reineke for his support and for very helpful discussions. 
\section{Notation and terminology}
\subsection{Representations of quivers} 
Let $k$ be an algebraically closed field.
\begin{defi}A quiver $Q$ consists of a set of vertices $Q_0$ and a set of
arrows $Q_1$ denoted by $\alpha: i \rightarrow j$ for $i, j\in Q_0$. In this situation $i$ is called the tail and $j$ the head of the arrow $\alpha$.
A quiver is finite if $Q_0$ and $Q_1$ are finite.\\
An oriented cycle in $Q$ is a set of arrows
$\alpha_1:i_1\rightarrow i_2,\alpha_2:i_2\rightarrow i_3\ldots,\alpha_n:i_n\rightarrow i_{n+1}$ such that $i_1=i_{n+1}$.\\
\end{defi}
In the whole paper we only consider quivers without oriented cycles.\\
Define the abelian group
\[\mathbb{Z}Q_0=\bigoplus_{i\in Q_0}\mathbb{Z}i\]
and its monoid of dimension vectors $\mathbb{N}Q_0$.
We introduce a non-symmetric bilinear form called the Euler form on $\mathbb{Z}Q_0$.
Define
\[\langle d,e\rangle :=\sum_{i\in Q_0}d_ie_i-\sum_{\alpha:i\rightarrow j}d_ie_j.\] A finite-dimensional $k$-representation of $Q$ is given by a tuple
\[X = ((X_i)_{i\in Q_0} , (X_{\alpha})_{\alpha\in Q_1}: X_i\rightarrow X_j)\]
of finite-dimensional $k$-vector spaces and $k$-linear maps between them. The
dimension vector $\underline{\dim}X \in\mathbb{N}Q_0$ of $X$ is defined by
\[\underline{\dim}X =\sum_{i\in Q_0}\dim_k X_ii.\]
Let $d\in\mathbb{N}Q_0$ be a dimension vector. The variety $R_d(Q)$ of $k$-representations
of $Q$ with dimension vector $d$ is defined as the affine $k$-space
\[R_d(Q) =\bigoplus_{\alpha:i\rightarrow j}\mathrm{Hom}_k(k^{d_i},k^{d_j}).\]
The algebraic group
\[G_d =\prod_{i\in Q_0}Gl_{d_i}(k)\]
acts on $R_d(Q)$ via simultaneous base change, i.e.
\[(g_i)_{i\in Q_0}\ast (X_{\alpha})_{\alpha\in Q_1} = (g_jX_{\alpha}g_i^{-1})_{\alpha:i\rightarrow j}.\]
The orbits are in bijection with the isomorphism classes of $k$-representations
of $Q$ with dimension vector $d$.\\\\
In the space of $\mathbb{Z}$-linear functions $\mathrm{Hom}_{\mathbb{Z}}(\mathbb{Z}Q_0,\mathbb{Z})$ we consider the basis given
by the elements $i^{\ast}$ for $i\in Q_0$, i.e. $i^{\ast}(j) = \delta_{i,j}$ for $j\in Q_0$. Define
\[\mathrm{dim}:=\sum_{i\in Q_0}i^{\ast}.\]
After choosing $\Theta\in \mathrm{Hom}_{\mathbb{Z}}(\mathbb{Z}Q_0,\mathbb{Z})$, we define the slope function $\mu : \mathbb{N}Q_0\rightarrow\mathbb{Q}$ via
\[\mu(d) =\frac{\Theta(d)}{\dim(d)}.\]
The slope $\mu(\underline{\dim}X)$ of a representation $X$ of $Q$ is abbreviated to $\mu(X)$.
\begin{defi} A representation $X$ of $Q$ is semistable (resp. stable) if for
all proper subrepresentations $0\neq U \subsetneq X$ the following holds:
\[\mu(U)\leq \mu(X)\text{ (resp. }\mu(U) < \mu(X)).\]
\end{defi}
Denote by $R^{ss}_d(Q)$ the set of semistable points and by $R^s_d(Q)$ the set of stable
points in $R_d(Q)$. In this situation we have the following theorem based on
Mumford's GIT and proved by King, see \cite{mum} and \cite{king}:
\begin{satz} We have:
\begin{enumerate}
\item The set of stable points $R^s_d(Q)$ is an open subset of the set of semistable
points $R^{ss}_d(Q)$, which is again an open subset of $R_d(Q)$.
\item There exists a categorical quotient $M^{ss}_d(Q):=R^{ss}_d(Q)//G_d$. Moreover,
$M^{ss}_d(Q)$ is a projective variety.
\item There exists a geometric quotient $M^s_d(Q):=R^s_d(Q)/G_d$, which is an
open smooth subvariety of $M^{ss}_d(Q)$.
\end{enumerate}
\end{satz}
\begin{bem}\label{dimension}
\end{bem}\begin{itemize}
\item For a stable representation $X$ its orbit in $R_d(Q)$ is of maximal possible
dimension, see \cite{king}. Since the scalar matrices act trivially on $R_d(Q)$,
the isotropy group is at least of dimension one. Therefore, if the moduli space is not empty, we get for its dimension
\[ \dim M^s_d(Q) = 1 - \langle d, d\rangle.\]
\end{itemize}
\subsection{Moduli spaces of vector bundles on the projective plane}
In this section we treat basic results concerning the theory of vector bundles
on projective spaces based on $\cite{oss}$, \cite{lep} and \cite{har}.\\
Let $\mathcal{E}\neq 0$ be a vector bundle on $\mathbb{P}^n$. Denote by $c_i(\mathcal{E})$ its $i$-th Chern class
and by $\mathrm{rk}(\mathcal{E})$ its rank. The Chern polynomial of $\mathcal{E}$ is given by \[c(\mathcal{E}) =\sum_{i=0}^{\mathrm{rk}(\mathcal{E})}c_i(\mathcal{E})t^i.\]
Define by
\[\mu(\mathcal{E}) = \frac{c_1(\mathcal{E})}{\mathrm{rk}(\mathcal{E})}\]
the slope of the vector bundle.
\begin{defi} A vector bundle $\mathcal{E}$ on $\mathbb{P}^n$ is semistable if we have
\[\mu(\mathcal{F}) \leq\mu(\mathcal{E})\]
for all coherent subsheaves $0 \neq \mathcal{F}\subseteq \mathcal{E}$.\\
The bundle $\mathcal{E}$ is stable if we have
\[\mu(\mathcal{F}) < \mu(\mathcal{E})\]
for all proper coherent subsheaves $\mathcal{F}$ with $0 < \mathrm{rk}(\mathcal{F}) < \mathrm{rk}(\mathcal{E})$.
\end{defi}
A stable bundle with respect to this stability condition is also often called
stable in the sense of Mumford-Takemoto or simply $\mu$-stable.
\begin{defi} The discriminant of a stable bundle of rank $r$ on $\mathbb{P}^2$ with
Chern classes $c_1$ and $c_2$ is defined by
\[D=\frac{1}{2r^2}(2rc_2-(r-1)c^2_1).\]
\end{defi}
In this situation we have the following theorem, see $\cite{bog}$:
\begin{satz} Let $\mathcal{E}$ be a semistable bundle of rank $r > 0$ and discriminant
$D$ on the projective plane. Then we have $D\geq 0$.
\end{satz}
There are lots of articles dealing with moduli spaces of stable vector bundles (or sheaves resp.) on projective varieties, see for instance \cite{sim}, \cite{mar1} or \cite{mar2}. In this paper we are only interested in moduli spaces of $\mu$-stable bundles on the projective plane. In particular, we are not interested in moduli spaces of semi-stable bundles. We do not describe the construction in detail and only need the following property of moduli spaces of $\mu$-stable vector bundles on the projective plane, again see for instance \cite{sim}, \cite{mar1} or \cite{mar2}:
\begin{satz}Fix the first two Chern classes $c_1$, $c_2$ and the rank $r$. Then
there exists a quasi-projective variety $\mathcal{M}(r, c_1, c_2)$ which parametrizes the isomorphism classes of $\mu$-stable vector bundles $\mathcal{E}$ on the projective plane such that $c_1(\mathcal{E})=c_1$, $c_2(\mathcal{E}) = c_2$ and $\mathrm{rk}(\mathcal{E}) = r$.
\end{satz}
In the following denote by $\mathcal{O}(1)$ the hyperplane bundle on the projective
plane and its dual line bundle by $\mathcal{O}(-1)$. Also define
\[\mathcal{O}(k) = \mathcal{O}(1)^{\otimes k}\]
if $k\geq 0$ and
\[\mathcal{O}(k) = \mathcal{O}(-1)^{\otimes k}\]
if $k<0$.\\
If we consider a bundle of rank $r$ on $\mathbb{P}^n$, its Chern classes $c_t$ with $t\in\mathbb{N}$ satisfy $c_t = 0$ for $t > \min(r, n)$. For details see for instance $\cite{lep}$.\\
Let $\mathcal{E}$ and $\mathcal{F}$ be vector bundles of rank $r$ and $s$ respectively on the projective space $\mathbb{P}^n$ and let
\[c(\mathcal{E}) =\prod_{i=1}^{r}(1+a_it)\]
and
\[c(\mathcal{F}) =\prod_{i=1}^{s}(1+b_it)\]
respectively with $a_i, b_i\in\mathbb{Z}$ be their Chern polynomials. Define $c_i$ with
$i = 0,\ldots, r + s$ by
\[\prod_{i,j}(1 + (a_i + b_j)t) =\sum_{i=0}^{r+s}c_it^i.\]
For the Chern polynomial of the bundle $\mathcal{E}\otimes\mathcal{F}$ we have, see \cite{hir} or \cite{har}:
\[c(\mathcal{E}\otimes\mathcal{F})=\sum_{i=0}^{\min(n,\mathrm{rk}(\mathcal{E}\otimes\mathcal{F}))}c_it^i.\]
If \[c(\mathcal{E}) =\sum_{i=0}^{\mathrm{rk}(\mathcal{E})}c_i(\mathcal{E})t^i,\] for the Chern polynomial of the dual bundle $\mathcal{E}^{\ast}$ we have
\[c(\mathcal{E}^{\ast})=\sum_{i=0}^{\mathrm{rk}(\mathcal{E})}c_i(\mathcal{E})(-t)^i.\]
Thus we obtain the following important property for the moduli spaces of
stable bundles on the projective plane:
\begin{pro}\label{pro1}
Twisting with a line bundle $\mathcal{O}(k)$ induces an isomorphism
of moduli spaces
\[\mathcal{M}(r, c_1, c_2)\cong \mathcal{M}(r, c_1 + rk, c_2 + (r - 1)kc_1 + k^2\frac{r(r - 1)}{2}).\]
If $r\leq 3$, the moduli spaces only depend on the discriminant.
\end{pro}
{\it Proof.} Let $\mathcal{E}$ be a stable bundle with Chern classes $c_1$, $c_2$ and rank $r$. The rank of $\mathcal{E}$ does not change after tensoring by $\mathcal{O}(k)$. From
\begin{eqnarray*}
c(\mathcal{E}\otimes\mathcal{O}(k))& = &(1 + (k + \frac{c_1}{2}-\sqrt{\frac{c_1^2}{4}-c_2})t)\\
&&(1 + (k + \frac{c_1}{2}+\sqrt{\frac{c_1^2}{4}-c_2})t)(1+kt)^{r-2}\modu t^3
\end{eqnarray*}
it follows
\[c(\mathcal{E}\otimes\mathcal{O}(k)) = 1+(rk + c_1)t+(\frac{r(r-1)}{2}k^2 +(r-1)kc_1+c_2)t^2.\]
Now an easy calculation implies that the discriminants of $\mathcal{E}$ and $\mathcal{E}\otimes\mathcal{O}(k)$ coincide. Furthermore, it is easy to see that (semi-)stable bundles remain (semi-)stable when tensoring with a line bundle. For $r=1,2$ this suffices to prove the second statement.\\ 
Thus let $r=3$. Because of the first part we may assume that $0\leq c_1\leq 2$. If $c_1=2$ we may take the dual and tensor with $\mathcal{O}(1)$ afterwards in order to get $c_1=1$. Note that the discriminant does not change when taking the dual. Under the assumption $c_1\in\{0,1\}$ it is easy to see that the moduli spaces only depend on the discriminant. 
\qed
Fixing the rank to be $r$, these moduli spaces will be denoted by $\mathcal{M}(r,D)$ or simply
$\mathcal{M}(D)$ if $r\leq 3$. By considering a torus action on the moduli space we can calculate
the Euler characteristic via reduction to torus fixed points.\\
The $(n + 1)$-dimensional torus $T = (\mathbb{C}^{\ast})^{n+1}$ acts on the projective space $\mathbb{P}^n$ via multiplication, i.e.
\[t\cdot(x_0 : x_1 :\ldots : x_n) = (t_0x_0 : t_1x_1 :\ldots : t_nx_n)\]
for $t = (t_0, t_1,\ldots, t_n)\in T$.
\begin{defi} A vector bundle $p : \mathcal{E}\rightarrow \mathbb{P}^n$ is toric if $T$ acts linear on the fibres such that for each $t\in T$ the following diagram commutes:
\[
\begin{xy}
\xymatrix@R30pt@C40pt{\mathcal{E}\ar[r]^t\ar[d]^p&\mathcal{E}\ar[d]^p\\
\mathbb{P}^n\ar[r]^t & \mathbb{P}^n
}\end{xy}
\]
\end{defi}
Let $E$ be a vector space. Then a descending $\mathbb{Z}$-filtration is defined as a chain
of subspaces $E(i)\subset E$, $i\in\mathbb{Z}$, such that $E(i)\subset E(i - 1)$. Denote the set
of all filtrations of a vector space $E$ by $\mathcal{F}(E)$. Analogously, we can define
families of filtrations of a space $E$, i.e. $E^{\alpha}\in\mathcal{F}(E)$ for $\alpha\in I$ and an index set $I$.\\\\
Let $\mathcal{E}$ be a toric bundle on $\mathbb{P}^2$. The $T =(\mathbb{C}^{\ast})^3$-action has an open orbit containing all points $p = (x_{\alpha} : x_{\beta} : x_{\gamma})\in\mathbb{P}^2$ such that $x_{\alpha},x_{\beta},x_{\gamma}\neq 0$. Let $E:= \mathcal{E}(p_0)$ be the fibre of an arbitrary point $p_0$ in this orbit. Since $\mathcal{E}$ is toric, we have $te\in\mathcal{E}(tp_0)$ for all $e\in E$. Now choose a generic point $p_{\alpha}$ from
the coordinate line
\[X_{\alpha} = \{(x_{\alpha} : x_{\beta} : x_{\gamma})\in\mathbb{P}^2\mid x_{\alpha} = 0\}.\]
Define
\[E^{\alpha}(i):= \{e \in E\mid\lim_{tp_0\rightarrow p_{\alpha}}\left(\frac{t_{\alpha}}{t_{\beta}}\right)^{-i}(te)\text{ exists}\}\]
and define $E^{\beta}, E^{\gamma}$ analogously. Obviously this definition is independent of
the choices of $p_0$ and $p_{\alpha}$. Instead of the chosen rational function, we may
also consider every other function with a pole of order $i$ in $X_{\alpha}$.\\
By this procedure for every bundle we get a family of descending $\mathbb{Z}$-filtrations
of $E$. Indeed we obviously have:
\[\ldots\subset E^{\alpha}(i + 1) \subset E^{\alpha}(i)\subset E^{\alpha}(i - 1)\subset\ldots\]
with the additional property $E^{\alpha}(i) = 0$ for $i\gg 0$ and $E^{\alpha}(i)=E$ for $i\ll 0$. From $\cite{kly2}$ we get the following theorem:
\begin{satz}\label{klyachko}
The category of toric bundles on the projective plane is equivalent to the category of vector spaces with a family of descending $\mathbb{Z}$-filtrations $E^{\alpha}\in\mathcal{F}(E)$ with $\alpha=1,2,3$ such that
\[ E^{\alpha}(i)=0\text{ for } i\gg 0\text{ and } E^{\alpha}(i) = E\text{ for } i\ll 0.\]
\end{satz}
Two filtrations $E_1,E_2\in\mathcal{F}(E)$ are isomorphic if there exists a $g\in Gl(E)$
such that $gE_2(i) = E_1(i)$ for all $i\in\mathbb{Z}$. Similarly, two families $E^{\alpha}_1, E^{\alpha}_2$ with $\alpha\in I$ for some index set $I$ are isomorphic if there exists a $g\in Gl(E)$ such that $gE^{\alpha}_1(i) = E^{\alpha}_2(i)$ for each $\alpha\in I$ and each $i\in\mathbb{Z}$.\\\\
Let $\mathcal{E}$ be a toric bundle on $\mathbb{P}^2$ given as a filtration $E^{\alpha}$, $\alpha = 1, 2, 3$. The first two Chern classes of this filtration and the corresponding bundle respectively are given as follows, see $\cite{kly2}$:
\[ c_1(\mathcal{E})=\sum_{i\in\mathbb{Z},\alpha}i\dim E^{[\alpha]}(i)\]
where $E^{[\alpha]}(i) = E^{\alpha}(i)/E^{\alpha}(i + 1)$ and
\[ c_2(\mathcal{E}) = \frac{c_1(\mathcal{E})^2}{2}-\frac{1}{2}\sum_{i\in\mathbb{Z},\alpha}i^2\dim E^{[\alpha]}(i)
-\sum_{\alpha\neq\beta,(i,j)\in\mathbb{Z}^2}ij\dim E^{[\alpha\beta]}(i,j)\]
where $E^{[\alpha\beta]}(i,j)=E^{\alpha}(i)\cap E^{\beta}(j)/(E^{\alpha}(i+1)\cap E^{\beta}(j)+ E^{\alpha}(i) \cap E^{\beta}(j + 1))$. The twist by again a line bundle corresponds to a shift of indices on the level of filtrations. More detailed we have the following:
\begin{lem}\label{lem1}
Let $\dim E = r$ and $k_{\alpha}\in\mathbb{Z}$ for $\alpha = 1,2,3$. By a shift of
indices $f:\mathbb{Z}^3\rightarrow\mathbb{Z}^3$, $(i_{\alpha}\rightarrow i_{\alpha}+k_{\alpha})_{\alpha=1,2,3}$, the discriminant remains constant. In particular, we may assume that the filtrations are in standard position, i.e. for $\alpha=1,2,3$ we have $E^{\alpha}(i)=E$ for all $i\geq 0$ and $E^{\alpha}(i)\neq E$ for all $i < 0$.
\end{lem}
{\it Proof.} Let $D'$ be the discriminant resulting from the index shift. Since
\[\sum_{i}\dim E^{[\alpha]}(i)=r,\]
it follows
\begin{eqnarray*}
D' &=& D +\sum_{\alpha}k^2_{\alpha}r^2+2(\sum_{\alpha}k_{\alpha}r)(\sum_{\alpha,i}i\dim E^{[\alpha]}(i))\\
&&+2\sum_{\alpha\neq\beta}k_{\alpha}k_{\beta}r^2-r\sum_{\alpha,i}2k_{\alpha}i\dim E^{[\alpha]}(i) -\sum_{\alpha}r^2k^2_{\alpha}\\
&&-2r\sum_{\alpha\neq\beta,i,j}k_{\alpha}j\dim E^{[\alpha,\beta]}(i,j)+k_{\beta}i\dim E^{[\alpha,\beta]}(i,j)\\
&&+k_{\alpha}k_{\beta}\dim E^{[\alpha,\beta]}(i,j)\\
&=& D + 2r\sum_{\alpha}k_{\alpha}\sum_{\beta\neq\alpha,i}i\dim E^{[\beta]}(i)+2r^2\sum_{\alpha\neq\beta}k_{\alpha}k_{\beta}\\
&&-2r\sum_{\alpha\neq\beta,i,j}k_{\alpha}j\dim E^{[\alpha,\beta]}(i,j) + k_{\beta}i\dim E^{[\alpha,\beta]}(i,j)\\
&&+k_{\alpha}k_{\beta}\dim E^{[\alpha,\beta]}(i,j).
\end{eqnarray*}
Now choose $n$ big enough such that $E^{\alpha}(n) = E$ and $E^{\alpha}(-n)=0$ for all $\alpha$.
Let $d_{\alpha,\beta}(i,j)=\dim E^{\alpha}(i)\cap E^{\beta}(j)$. Obviously we have
\[\dim E^{[\alpha\beta]}(i,j)=d_{\alpha,\beta}(i,j)-d_{\alpha,\beta}(i+1,j)-d_{\alpha,\beta}(i,j+1)+d_{\alpha,\beta}(i+1,j+1).\]
For $\alpha\neq\beta$ we have
\begin{eqnarray*}
\sum_{i,j}\dim E^{[\alpha,\beta]}(i,j) &=&\sum_{i=-n}^{n}\sum_{j=-n}^n \dim E^{[\alpha,\beta]}(i,j)\\
&=& d_{\alpha,\beta}(-n,-n)-d_{\alpha,\beta}(-n,n+1)\\
&&-d_{\alpha,\beta}(n+1,-n) + d_{\alpha,\beta}(n+1,n+1)\\
&=& r.
\end{eqnarray*}
Then we have $d_{\alpha,\beta}(i,j)=0$ if $i \gg 0$ or $j\gg 0$. Moreover, it follows
\[\sum_jj\dim E^{[\beta]}(j) =\sum_{j=-n}^n j(d_{\alpha,\beta}(-n,j) - d_{\alpha,\beta}(-n,j+1)\]
and we get analogously to the previous equation
\begin{eqnarray*}
\sum_{i,j}j\dim E^{[\alpha,\beta]}(i,j)&=&\sum_{i=-n}^n\sum_{j=-n}^nj(d_{\alpha,\beta}(i,j)-d_{\alpha,\beta}(i+1,j)\\
&&-d_{\alpha,\beta}(i,j+1)+d_{\alpha,\beta}(i+1,j+1))\\
&=&\sum_{j=-n}^nj(d_{\alpha,\beta}(-n,j)-d_{\alpha,\beta}(-n,j+1)\\
&&-d_{\alpha,\beta}(n+1,j)+d_{\alpha,\beta}(n+1,j+1)\\
&=&\sum_{j=-n}^n j(d_{\alpha,\beta}(-n,j) - d_{\alpha,\beta}(-n,j+1)).
\end{eqnarray*}
\qed
We consider the subspace quiver with the vertex set
\[Q_0 =\{q_0\}\cup\{q_{i,j}\mid 1\leq i\leq n,j\in\mathbb{N}^+\}\]
and arrow set
\[Q_1 =\{\alpha: q_{i,1}\rightarrow q_0 \mid 1 \leq i\leq n\}\cup\{\alpha:q_{i,j+1}\rightarrow q_{i,j}\mid 1\leq i\leq n,j\in\mathbb{N}^+\}.\]
Let $X$ be a representation with dimension vector $d=(d_i)_{i\in Q_0}$ such that
$d_{i,j+1}\leq d_{i,j}$. Denote the linear maps corresponding to $X$ by $X_{i,j}:\mathbb{C}^{d_{i,j+1}}\rightarrow\mathbb{C}^{d_{i,j}}$, where we assume that all maps $X_{i,j}$ are injective. It is easy to see that every such representation is isomorphic to a representation $X'$ such that
\[X'_{i,j}=\binom{E_{d_{i,j}}}{0}\]
for all $j\geq 1$, where $E_{d_{i,j}}$ is the $d_{i,j}\times d_{i,j}$-identity matrix. In what follows we assume that all representation are of this type.\\
Thus a representation $X$ is given by a $n$-tuple of matrices $(X_{i,0})_{1\leq i\leq n}$ and a
dimension vector $(d_i)_{i\in Q_0}$. Given a representation $X$ we get a filtration as follows: let $(X_{i,0})_k$ be $k$-th column and define
\[E^i(j) = \langle (X_{i,0})_{n-d_{i,j}+1},\ldots, (X_{i,0})_n\rangle.\]
If two filtration are isomorphic, the corresponding two representations are obviously isomorphic as well (via the same $g$).\\
The other way around, we obtain a representation from a filtration. If two representations $X,X'$ are isomorphic, there exists a $g\in Gl_d =\Pi_{i\in Q_0}Gl_{di}$
with $g\ast X = X'$. If $d_{i,j} < d_{i,j-1}$, we have
\[g_{i,j-1}=\begin{pmatrix}
g_{i,j}&\star\\
0&\star\end{pmatrix}\]
for all $j\geq 2$ where $g_{i,j}\in\mathbb{C}^{d_{i,j}\times d_{i,j}}$. In particular, there exists a matrix $g_0\in Gl_{d_0}(E)$ and matrices $g_{i,1}$ such that
\[g_0X_{i,1} = X'_{i,1}g_{i,1}\]
for all $1 \leq i\leq n$ where the $g_{i,1}$ are nested such that the subspaces $E^i(j)$ of
the corresponding filtration are invariant under $g_{i,j}$.\\\\
\begin{bem}
\end{bem}
\begin{itemize}
\item If we in general consider vector bundles on the projective space $\mathbb{P}^n$
such that $n\geq 3$, Klyachko's theorem $\ref{klyachko}$ just holds on an additional
condition. The arms are in bijection with vectors generating the fan belonging to $\mathbb{P}^n$. Thereby we consider the projective space as a toric variety. In these cases we get a toric bundle from a filtration if all subfiltrations belonging to those arms, which correspond to vectors generating a cone of the fan, generate a distributive lattice. This is
automatically satisfied in the case $n = 2$. If $n\geq 3$ this already means that every toric bundle of rank $2$ splits.
\end{itemize}
We consider the stability condition given by the slope function
\[\mu(d) =\frac{\Theta(d)}{\dim d}\]
with $\Theta = -q_0^{\ast}$. The following holds, see $\cite{kly1}$:
\begin{satz}
Let $\mathcal{E}$ be a toric bundle on $\mathbb{P}^2$ given by a triple of filtration $E^{\alpha}$. Then the following are equivalent:
\begin{enumerate}
\item $\mathcal{E}$ is stable in the sense of Mumford-Takemoto.
\item The family of subspaces $E^{\alpha}(i) \subset E$ is stable under the action of $Gl(E)$
in the sense of Mumford.
\item For all subspaces $0 \subset F \subset E$ we have
\[\sum_{\alpha,i>N}\frac{\dim E^{\alpha}(i)\cap F}{\dim F}<\sum_{\alpha,i>N}\frac{\dim E^{\alpha}(i)}{\dim E}.\]
\end{enumerate}
\end{satz}
{\it Proof.} It is easy to see that the introduced stability for the subspace quiver
is equivalent to the third assertion. Therefore, we have the equivalence with
two. Furthermore, it follows from the definition of the stability via first
Chern class and rank that the first and third statements are equivalent.
\qed
Denote by $\mathcal{U}(r)$ the set of all subspace quivers with three arms, i.e. $n = 3$,
with dimension vector $d$ such that $d_{q_0} = r$ and $d_{q_{i,j}}\geq d_{q_{i,j+1}}$.\\ 
The preceding theorem means that the moduli spaces of stable representations and the moduli spaces of stable filtrations of fixed length are isomorphic. The former will be investigated in greater detail.\\ In particular, if $r\leq 3$ every point in $\mathcal{M}(D)^T$ corresponds to a stable representation (up to isomorphism) of the subspace quiver $\mathcal{U}(r)$, see $\ref{pro1}$ and $\ref{lem1}$. Moreover, we may understand fixed point components of $\mathcal{M}(D)$ as moduli spaces of the subspace quiver. \\Note that for $r=2$ every point in $\mathcal{M}(D)^T$ uniquely corresponds to a stable representation (up to isomorphism) of $\mathcal{U}(2)$.\\
For $r=3$ and $D\equiv 0\modu 6$, this correspondence is unique as well. Indeed, such a representation (filtration in standard position resp.) can be shifted such that we get a filtration with $c_1=0$.\\
For $r=3$, $D\equiv 4\modu 6$ and $c_1\equiv 1\modu 3$, such a representation can be shifted such that we get a filtration with $c_1=1$. Analogously, for $c_1\equiv 2\modu 3$ we get a corresponding filtration with $c_1=2$. Because of $\ref{pro1}$ this means that every point in $\mathcal{M}(D)^T$ corresponds to exactly two stable representation of $\mathcal{U}(3)$. 
\section{Systems of linear inequalities and polyhedrons}\label{LU}
In this section a summary of required methods concerning systems of linear
inequalities is given. For more details see for instance $\cite{pad}$, which also
provides the basis of this section. We will not place emphasise on finding a
solution to the system of linear inequalities as efficient as possible because
our applications in the next sections do not need it. Thus we just discuss
how to get a solution of a given system of linear inequalities.
\begin{defi}
Let $A\in\mathbb{R}^{m,n}$ be a $(m\times n)$-matrix and $b\in\mathbb{R}^n$. A polyhedron
$P\subseteq\mathbb{R}^n$ is the set of solutions $x\in\mathbb{R}^n$ of some system of linear inequalities $Ax\leq b$.
\end{defi}
In the following we denote the polyhedron coming from a matrix $A$ and a vector $b$ by $P(A,b)$.\\
A polytope is a bounded polyhedron. This means there exists some $s\in\mathbb{R}$ such that for every $x\in P(A,b)$ we have $\parallel x \parallel\leq s$.\\
Thus a polyhedron is determined by the solution of a finite number of inequalities. Every inequality defines some half space so that a polyhedron can be understood as the intersection of a finite number of half spaces.
\begin{defi} Let $x_1, x_2,\ldots,x_k\in\mathbb{R}^n$. The convex hull of these points is
defined by
\[\mathrm{conv}(x_1,\ldots,x_k):=\{\sum_{i=1}^k\mu_ix_i\mid\sum_{i=1}^k\mu_i=1,\mu_i\geq 0\}.\]
The convex cone is defined by
\[\mathrm{cone}(x_1,\ldots,x_k):=\{\sum_{i=1}^k\mu_ix_i\mid\mu_i\geq 0\}.\]
\end{defi}
\begin{defi}
Let $P\subset\mathbb{R}^n$ be some subset. A point $x\in P$ is called an
extreme point of $P$ if for all $x_1, x_2\in P$ and every $0 <\mu <1$ such that
$x = \mu x_1 + (1-\mu)x_2$, we have $x = x_1 = x_2$.
\end{defi}
Thus $x$ can be uniquely written as a convex combination of elements of $P$,
namely as the trivial one.\\
The following theorem plays an important role if we want to determine the
solutions of a given system of linear inequalities, for a proof see $\cite{pad}$.
\begin{satz}
Let $A\in \mathbb{R}^{m,n}$ be a $(m\times n)$-matrix and $b\in\mathbb{R}^n$ such that
$m\geq n$. A point $x_0\in P(A,b)$ is an extreme point of the polyhedron if
$Ax_0\leq b$ and $A'x_0=b'$ for some $(n\times n)$-submatrix of $A$ with $\mathrm{rank}(A') = n$
and the corresponding subvector $b'$ of $b$.
\end{satz}
By corresponding subvector we mean of course that $b'$ results from $b$ as follows: we remove the entry $b_i$ if and only if we remove the $i$-th row of $A$.\\
In addition to the extreme points, whose convex hull corresponds to a polytope satisfying the inequalities, we determine some vectors whose positive linear combinations based on this polytope describe all solutions of the given inequalities.\\
Thus let $x_0,x\in\mathbb{R}^n$. Consider $x_{\lambda} = x_0 + \lambda x$ for $\lambda\geq 0$. Then we have $x_{\lambda}\in P(A,b)$ for every $\lambda\geq 0$ if and only if $x_0\in P(A,b)$ and $Ax\leq 0$.\\
This leads us to the following definitions:
\begin{defi}
\begin{enumerate}
\item A set $C\subseteq\mathbb{R}^n$ is a cone if for every pair of points $x_1, x_2\in C$ we have $\lambda_1x_1 + \lambda_2x_2\in C$ for all $\lambda_1,\lambda_2\geq 0$.
\item A cone is called pointed if it does not contain any subspace except $\{0\}$.
\item A half-line $y = \{\lambda x\mid \lambda\geq 0,x\in\mathbb{R}^n\}$ is an extremal ray of $C$ if $y\in C$ and $-y\notin C$ and if for all $y_1,y_2\in C$ and $0 <\mu < 1$ with
$y = (1 -\mu)y_1 + \mu y_2$ we already have $y = y_1 = y_2$.
\end{enumerate}
\end{defi}
The polyhedral cone corresponding to some system of inequalities $(A,b)$ is
defined as
\[C(A) = \{x\in\mathbb{R}^n\mid Ax\leq 0\}.\]
Obviously $C(A)$ is both a polyhedron and a cone.\\
The following theorem describes how to determine all extremal rays of some
polyhedral cone. For a proof again see \cite{pad}.
\begin{satz}
Let $C(A)$ be a pointed cone. Then $x\in C(A)$ is an extremal
ray of $C(A)$ if and only if there exist $\mathrm{rank}(A) - 1$ linear independent row
vectors $a_1, \ldots,a_{\mathrm{rank}(A)-1}$ of $A$ such that
\[\left(\begin{array}{c}
a_1\\ \vdots\\ a_{\mathrm{rank}(A)-1}
\end{array}\right)\cdot x=0\]
and moreover $Ax\leq 0$ holds.
\end{satz}
In what follows we assume that every polyhedron does not contain any one-
dimensional subspace. Therefore, we only consider systems of inequalities
whose solutions $x\in\mathbb{R}^n$ satisfy the additional condition $x_i\geq 0$. Such a
system is said to be in standard form.\\
This assumption is no restriction because every system can be transformed
into a system in standard form. The advantage of such a system is that the
set of solutions does not contain lines, i.e. the corresponding polyhedron
and in particular the corresponding polyhedral cone are pointed.\\
Moreover, note that the sets of extreme points and extremal rays are finite,
which is clear because of the preceding theorems.\\
In conclusion we have the following:
\begin{satz}
Let $(A,b)$ be a system of inequalities in standard form. Let
$X =\{x_1, x_2,\ldots,x_s\}$ be the set of all extreme points of the polyhedron $P(A, b)$
and $Y = \{y_1,y_2,\ldots, y_t\}$ the set of all extremal rays of the polyhedral cone
$C(A)$. Then the polyhedron $P(A,b)$ consisting of all solutions of the system
of linear inequalities defined by $Ax\leq b$ is given by
\[P(A, b) = \mathrm{conv}(X) + \mathrm{cone}(Y).\]
\end{satz}
\section{Euler characteristic of moduli spaces of stable bundles}
\subsection{The case of rank two bundles on the projective plane}
In this section we first review the methods presented in $\cite{kly1}$ in order to
deduce from it a similar formula for the Euler characteristic of rank three
bundles on the projective plane. In order to compare both results we make
a small modification of Klyachko's methods.\\
Denote by $H(D)$ the Hurwitz function counting the number of classes of
reduced binary quadratic forms $Q$ with discriminant $D$ with weight $\frac{2}{|\mathrm{Aut} Q|}$. Referring to \cite{kly1} the Euler characteristic of moduli spaces of stable bundles
of rank two on the projective plane is given by
\[\chi(\mathcal{M}(c_1,c_2))=\left\{\begin{array}{l} 3H(D),\text{ if } D\equiv -1\modu 4\\3H(D)-\frac{3}{2}d(\frac{D}{4}),\text{ if } D\equiv 0 \modu 4
\end{array}.\right .\]
The starting point for the derivation of this formula is the following well-known theorem, see for instance \cite{car}, \cite{cg} or \cite{em}:
\begin{satz} Let $X$ be a complex variety on which a torus $T$ acts. For the Euler characteristic $\chi$ of $X$ we have
\[\chi(X) = \chi(X^T).\]
\end{satz}
Denote by $E^k$ with $k = 1, 2, 3$ a triple of filtrations of some two-dimensional vector space $E$,
i.e. $E^k$ consists of $E$ and filtrations $E^k(i)$ with $k = 1, 2, 3$ and $i\in\mathbb{Z}$ with
the additional condition
\[E^k(i) = 0\text{ for } i \gg 0\text{ and }E^k(i) = E\text{ for }i \ll 0.\]
Fixing such a triple, we define
\[\alpha_k := \mid\{i\mid\dim E^k(i) = 1\}\mid\]
for $k = 1, 2, 3$. The stability condition corresponds to the inequalities
\begin{equation}\label{ungl1}\alpha_1<\alpha_2+\alpha_3\text{ and }\alpha_2<\alpha_1+\alpha_3\text{ and }
\alpha_3<\alpha_1+\alpha_2.\end{equation}
Therefore, for every triple $(\alpha_1, \alpha_2, \alpha_3)$ satisfying these inequalities, there exists at least one stable bundle. Note that we always assume that the filtrations are in standard form.\\
Hence the discriminant $-D = c^2_1- 4c_2$ is given by
\[-D=\alpha_1^2+\alpha_2^2+\alpha_3^2-2\alpha_1\alpha_2-2\alpha_2\alpha_3-2\alpha_1\alpha_3.\]
As proved in the last section, every filtration corresponds to a representation
of a subspace quiver which is the quiver with three arms meeting in one
point in this case. Since there exists a stable representation, by use of the
dimension formula \ref{dimension} we get that the moduli spaces of such quivers are
zero-dimensional. Therefore, the Euler characteristic is one in each case.
Thus we have
\[\chi(\mathcal{M}(c_1,c_2)) = \mid\mathcal{M}(c_1, c_2)^T\mid.\]
By considering the inequalities (\ref{ungl1}) and applying the theorems of the first
section, we get that $P = (1, 1, 1)$ is the only extreme point of the system of
linear inequalities (\ref{ungl1}).\\
Moreover, we obtain the extremal rays $v_1 = (1, 1, 0)$, $v_2 =(0,1,1)$ and
$v_3 = (1, 0, 1)$. Therefore, we get that all positive integer valued solutions are
of the form
\begin{eqnarray*}v&=&(1,1,1)+k_1(1,1,0)+k_2(0,1,1)+k_3(1,0,1)\\&=&(k_1+k_3+1,k_1+k_2+1,k_2+k_3+1)\\
\end{eqnarray*} with $k_1,k_2,k_3\in\mathbb{Q}^+$.\\
Define
\[\mathbb{L}=\{v\in\mathbb{N}^3\mid v=(k_1+k_3+1,k_1+k_2+1,k_2+k_3+1),k_1,k_2,k_3\in\mathbb{Q}^+\}.\]
Let $k_1, k_2, k_3 < 1$ and $k_i \neq 0$ for at least one $k_i$. Obviously the only solution
we obtain in this way is $(2, 2, 2)$. Indeed, that is the case if $k_i =\frac{1}{2}$ for all
$i =1, 2, 3$. From this we get
\[\mathbb{L}=\{v\in\mathbb{N}^3\mid v=(k_1+k_3+i,k_1+k_2+i,k_2+k_3+i),k_1,k_2,k_3\in\mathbb{N},i=1,2\}.\]
Indeed, if $v \in\mathbb{L}$ is a solution with $k_i\in\mathbb{Q}+$ for $i = 1, 2, 3$, there also exists a solution $v' \in\mathbb{L}$ with $v'=(k'_1,k'_2,k'_3)$, where $k_i':=k_i-\lfloor k_i\rfloor$. Following the
consideration from above we have $v' =(1, 1, 1)$ or $v' = (2, 2, 2)$.
Now consider
\[\mathbb{L}_i=\{v\in\mathbb{N}^3\mid v=(k_1+k_3+i,k_1+k_2+i,k_2+k_3+i),k_1,k_2,k_3\in\mathbb{N}\}\]
with $i = 1, 2$. It is easy to see that $\mathbb{L}_1\cap\mathbb{L}_2 =\emptyset$. If $v\in\mathbb{L}_1$, for the discriminant we get
\[-D=-4k_1k_2-4k_1k_3-4k_2k_3-4k_1-4k_2-4k_3-3\]
and for $v\in\mathbb{L}_2$ we have
\[-D=-4k_1k_2-4k_1k_3-4k_2k_3-8k_1-8k_2-8k_3-12.\]
In particular, this means that every solution $v\in\mathbb{L}_1$ belongs to a moduli
space with $D\equiv 3\modu 4$ and every solution $v\in\mathbb{L}_2$ to a moduli space with
$D\equiv 0\modu 4$.\\
Obviously all solutions are uniquely determined by $k_1, k_2$ and $k_3$. Note that
the second equation is equivalent to the diophantine equation
\[xy + yz + zx = n\]
for $x, y, z\geq 1$ and $n\in\mathbb{N}$ with $n\geq 3$. We can see this by dividing by $-4$ and
defining $x = k_1 - 1, y = k_2-1$ and $z = k_3 - 1$ afterwards. For more details
concerning this diophantine equation see \cite{pet}.\\\\
In this case we get for the generating function of the Euler characteristic
\begin{eqnarray*}F(x)&=&\sum_{i=0}^\infty \chi(\mathcal{M}(-4i))x^{4i}\\
&=&\sum_{(k_1,k_2,k_3)\in\mathbb{N}_0^3}x^{4k_1k_2+4k_2k_3+4k_1k_3+8(k_1+k_2+k_3)+12}.\\
\end{eqnarray*}
In the other case we get
\begin{eqnarray*}F(x)&=&\sum_{i=0}^\infty \chi(\mathcal{M}(-4i-3))x^{4i+3}\\
&=&\sum_{(k_1,k_2,k_3)\in\mathbb{N}_0^3}x^{4k_1k_2+4k_2k_3+4k_1k_3+4(k_1+k_2+k_3)+3}.\\
\end{eqnarray*}
\subsection{The case of stable rank three bundles on the projective plane}
Let $\alpha_{ij}\geq 0$ with $i\in\{1, 2, 3\}$ and $j\in\{1,2\}$. We consider the subspace
quiver with dimension vectors defined by
\[\dim(q_0)=3,\]
\[\dim(q_{i,k})=2\text{ for } 1\leq i\leq 3,1\leq k\leq\alpha_{i2}\]
and
\[\dim(q_{i,k})=1\text{ for } 1\leq i\leq 3,\alpha_{i2}+1\leq k\leq\alpha_{i2}+\alpha_{i1}.\]
In the following denote this quiver by $\mathcal{U}(\alpha_{11},\alpha_{12},\alpha_{21},\alpha_{22},\alpha_{31},\alpha_{32})$.\\\\
In the following denote by $U_{ij}$ the six different subspaces with $i\in\{1,2,3\}$ and $j\in\{1,2\}.$ Obviously, we always have $U_{i1}\subset U_{i2}$. This means the first
Chern class is given as follows:
\[c_1(\mathcal{E})=\alpha_{11}+\alpha_{21}+\alpha_{31}+2\alpha_{12}+2\alpha_{22}+2\alpha_{32}.\]
By considering the second Chern class the following problem appears: if
we fix a quiver $\mathcal{U}(\alpha_{11},\alpha_{12},\alpha_{21},\alpha_{22},\alpha_{31},\alpha_{32})$, the second Chern class varies,
depending on the number of two-dimensional subspaces that contain the
one-dimensional subspaces.\\
We consider the cases $U_{i1}\nsubseteq U_{k2}$ for all $i$ and $k\neq i$ and $U_{i1}\subset U_{j2}$ for $i\neq j$.\\
We first assume $U_{i1}\nsubseteq U_{k2}$ for all $i$ and $k\neq i$. Then we have:
\[c_2(\mathcal{E})=\sum_{i=1}^3\alpha_{i2}^2+\alpha_{i1}\alpha_{i2}+\sum_{1\leq i<j\leq 3}\alpha_{i1}\alpha_{j1}+2\alpha_{i1}\alpha_{j2}+2\alpha_{i2}\alpha_{j1}+3\alpha_{i2}\alpha_{j2}.\]
In the following let $i,j$ and $k$ be mutually different. If we choose $U_{i1}$ and
$U_{i2}$ as subspaces in the stability condition with $1\leq i\leq 3$, we get the six
inequalities:
\[\alpha_{i1}+2\alpha_{i2}<2\alpha_{j1}+\alpha_{j2}+2\alpha_{k1}+\alpha_{k2}\]
and
\begin{equation}\label{eins}2\alpha_{i1}+\alpha_{i2}<\alpha_{j1}+2\alpha_{j2}+\alpha_{k1}+2\alpha_{k2}.\end{equation}
respectively. Considering the subspace $U_{i2}\cap U_{j2}$ for $i\neq j$ we get the following
condition:
\begin{equation}\label{zwei}\alpha_{i2}+\alpha_{j2}<\alpha_{11}+\alpha_{21}+\alpha_{31}+2\alpha_{k2}.\end{equation}
Choosing $U_{i1}\oplus U_{j1}$ we have
\begin{equation}\label{drei}\alpha_{i1}+\alpha_{j1}<\alpha_{12}+\alpha_{22}+\alpha_{32}+2\alpha_{k1}.\end{equation}
Obviously we do not have to consider other subspaces in order to test a
representation for stability. Thus in the case $U_{k1}\nsubseteq U_{l2}$ for each $l\neq k$, the
discriminant is given by
\begin{eqnarray*}D&=&2c_1^2-6c_2\\
&=&\sum_{i=1}^3 2\alpha_{i1}^2+2\alpha_{i1}\alpha_{i2}+2\alpha_{i2}^2\\
&&-2\left(\sum_{1\leq i<j\leq 3}\alpha_{i1}\alpha_{j1}+2\alpha_{i1}\alpha_{j2}+2\alpha_{i2}\alpha_{j1}+\alpha_{i2}\alpha_{j2}\right).\\
\end{eqnarray*}
Now consider the case $U_{i1}\subset U_{j2}$ for $i\neq j$. If $U_{i2} = U_{j2}$, we would get
\[\alpha_{i1}+2\alpha_{i2}+\alpha_{j1}+2\alpha_{j2}<2\alpha_{k1}+\alpha_{k2}\]
contradicting the inequalities (\ref{eins})-(\ref{drei}). If $U_{i1} = U_{j1}$, we analogously obtain the inequality
\[2\alpha_{i1}+\alpha_{i2}+2\alpha_{j1}+\alpha_{j2}<\alpha_{k1}+2\alpha_{k2}\]
again contradicting the above ones.\\
Thus it remains to consider the stability condition in the cases $U_{i1}\subset U_{j2}$.
Then we get the additional inequalities
\[2\alpha_{i1}+\alpha_{j2}+\alpha_{i2}<2\alpha_{k2}+\alpha_{k1}+\alpha_{i1}\]
and
\begin{equation}\label{vier}2\alpha_{j2}+\alpha_{i1}+\alpha_{j1}<2\alpha_{k1}+\alpha_{k2}+\alpha_{i2}.\end{equation}
Obviously, they do not conflict with the above inequalities. Actually, if we
have nowhere equality, we obtain that exactly four of these twelve inequalities have to be satisfied. We will shortly come back to this point.\\
First we consider the second Chern class in the case $U_{k1}\subseteq U_{l2}$ for $k\neq l$.
Then we obtain
\begin{eqnarray*}c_2(\mathcal{E})&=&\sum_{i=1}^3\alpha_{i2}^2+\alpha_{i1}\alpha_{i2}\\
&&+\left(\sum_{1\leq i<j\leq 3}\alpha_{i1}\alpha_{j1}+2\alpha_{i1}\alpha_{j2}+2\alpha_{i2}\alpha_{j1}+3\alpha_{i2}\alpha_{j2}\right)-\alpha_{k1}\alpha_{l2}.\end{eqnarray*}
Thus for the discriminant we obtain:
\begin{eqnarray*}D&=&2c_1^2-6c_2\\
&=&\sum_{i=1}^3 2\alpha_{i1}^2+2\alpha_{i1}\alpha_{i2}+2\alpha_{i2}^2\\
&&-2\left(\sum_{1\leq i<j\leq 3}\alpha_{i1}\alpha_{j1}+2\alpha_{i1}\alpha_{j2}+2\alpha_{i2}\alpha_{j1}+\alpha_{i2}\alpha_{j2}\right)+6\alpha_{k1}\alpha_{l2}.\\
\end{eqnarray*}
It is easy to see that in general the discriminant satisfies the property:
\[D\equiv 0\modu 6\text{ or } D\equiv 4\modu 6.\]
Also note that the second Chern class does not change if $U_{k1} \subset U_{i1} \oplus U_{j1}$
for mutually different $i$, $j$ and $k$.\\
Let $\alpha_{ij} = 1$ for all $i, j$. It is easy to see that every filtration is isomorphic to
one of the following form:
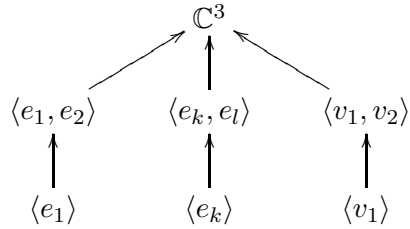
\filt{e_1}{e_2}{e_k}{e_l}{v_1}{v_2}
Thereby $k,l\in \{1,2,3\}$ such that $k\neq l$ and $v_1\neq v_2$ are arbitrary vectors.
Obviously we obtain the same for filtrations of arbitrary length. Such a
filtration is said to be in standard form.
\subsection{The case $D\equiv 4\modu 6$}
Let $\mathcal{U}(\alpha_{11},\alpha_{12},\alpha_{21},\alpha_{22},\alpha_{31},\alpha_{32})$ be a subspace quiver with $D\equiv 4\modu 6$.\\
Considering the inequalities $(\ref{eins})$-$(\ref{vier})$ we get the following lemma:
\begin{lem}
Let $\mathcal{U}(\alpha_{11},\alpha_{12},\alpha_{21},\alpha_{22},\alpha_{31},\alpha_{32})$ be the subspace quiver such that $D\equiv 4\modu 6$. Then we have
\begin{enumerate}
\item There exist no semistable points.
\item There exist exactly two stable points such that $U_{i1}\subseteq U_{j2}$ for $i\neq j$.
\end{enumerate}
\end{lem}
{\it Proof.} Let $\alpha\in\mathbb{N}^6$ such that there exists a semistable point for $\mathcal{U}(\alpha)$. We can without lose of generality assume that
\[\alpha_{11}+2\alpha_{12}=2\alpha_{21}+\alpha_{22}+2\alpha_{31}+\alpha_{32}.\]
For the discriminant we get in this case by a straight forward calculation
\begin{eqnarray*}D&=&6\alpha_{11}^2+6\alpha_{22}^2+6\alpha_{32}^2-6\alpha_{21}\alpha_{31}+6\alpha_{22}\alpha_{32}-6\alpha_{21}\alpha_{32}-6\alpha_{22}\alpha_{31}\\&&
-6\alpha_{11}\alpha_{22}-12\alpha_{11}\alpha_{32}.\end{eqnarray*}
This proves the first assertion.\\
The second part is proved as follows: if $U_{i1}\subseteq U_{j2}$, this point cannot be
semistable. Indeed, otherwise the discriminant belonging to this point would
satisfy $D\equiv 0\modu 6$ what is checked as before. By considering the twelve
inequalities (\ref{vier}) in detail we see that always exactly two pairs of them have
to be satisfied.
\qed
Let $M(\mathcal{U}(\alpha))^s$ the moduli space of stable representations of $\mathcal{U}(\alpha)$. It coincides with the moduli space of semistable representations. Following \cite{king} and \cite{kac}, therefore, it is a smooth projective variety of dimension one. By \cite{sch} we get that this projective curve (resp. the moduli space) is rational.
Thus it follows, see for instance \cite{har}, that
\[M(\mathcal{U}(\alpha))^s \cong\mathbb{P}^1.\]
Denote by $D(x)$ (resp. $c_1(x)$) the discriminant (resp. first Chern-class) corresponding to a point $x\in M(U(\alpha))^s$.
Moreover, define
\[M(\mathcal{U}(\alpha))^s_{D}=\{x\in M(\mathcal{U}(\alpha))^s\mid D(x)=D\}\]
and
\[M(\mathcal{U}(\alpha))^s_{D,i}=\{x\in M(\mathcal{U}(\alpha))^s\mid D(x)=D,\,c_1(x)\equiv i\modu 3\}.\]
By the preceding section and the preceding lemma, it follows that there exists exactly one tuple $(D_1,i)\in\mathbb{N}\times\{1,2\}$ such that $M(\mathcal{U}(\alpha))^s_{D_1,i}=\mathbb{P}^1\backslash\{(1:0),(0:1)\}$. Moreover, there exist exactly two tuples $(D_2,j),(D_3,k)\in\mathbb{N}\times\{1,2\}$ such that $M(\mathcal{U}(\alpha))^s_{D_i}=\{\mathrm{pt}\}$ for $i = 2, 3$. Then by the methods of the second chapter we obtain
\[M(D)^T\cong\bigcup_{\alpha\in\mathbb{N}^6}M(\mathcal{U}(\alpha))_{D,1}^s\cong\bigcup_{\alpha\in\mathbb{N}^6}M(\mathcal{U}(\alpha))_{D,2}^s.\]
Therefore, we get
\begin{eqnarray}\label{six}
\chi(M(D)^T)=\frac{1}{2}\sum_{\alpha\in\mathbb{N}^6}\chi(M(\mathcal{U}(\alpha))_D^s).\end{eqnarray}
Since $\chi(\mathbb{P}^1\backslash\{(1:0),(0:1)\})=0$, we just need to consider such moduli
spaces satisfying $M(\mathcal{U}(\alpha))^s_{D_i}=\{\mathrm{pt}\}$ with $i = 2, 3$. They correspond to the inclusions $U_{i1}\subseteq U_{j2}$.\\
If we consider the twelve inequalities (\ref{vier}) we get the following possibilities:
\begin{enumerate}
\item $U_{i1}\subset U_{j2},U_{j1}\subset U_{i2}$
\item $U_{i1}\subset U_{j2},U_{i1}\subset U_{k2}$
\item $U_{i1}\subset U_{j2},U_{k1}\subset U_{j2}$
\end{enumerate}
for mutually different $i$, $j$ and $k$. Therefore, the aim is to find the solutions of
the system of linear inequalities consisting of the inequalities $(\ref{eins})-(\ref{drei})$ and
the four inequalities $(\ref{vier})$ corresponding to these three cases.\\
Assume $\alpha_{ij}\neq 0$, the case $\alpha_{ij} = 0$ is discussed as a special case later.
\begin{defi}
Fix $\mathcal{U}(\alpha):=\mathcal{U}(\alpha_{11},\alpha_{21},\alpha_{31},\alpha_{12},\alpha_{22},\alpha_{32})$. Let
 \[\sigma_{12}\cdot (\alpha_{11},\alpha_{21},\alpha_{31},\alpha_{12},\alpha_{22},\alpha_{32}):=(\alpha_{21},\alpha_{11},\alpha_{31},\alpha_{22},\alpha_{12},\alpha_{32}).\]
Further define $\sigma_{13}, \sigma_{23}$ analogously and in addition
\[\tau\cdot (\alpha_{11},\alpha_{21},\alpha_{31},\alpha_{12},\alpha_{22},\alpha_{32}):=(\alpha_{12},\alpha_{22},\alpha_{32},\alpha_{11},\alpha_{21},\alpha_{31}).\]
\end{defi}
In doing so we get a group $G=\langle\sigma_{ij},1\leq i< j\leq 3,\tau\rangle$ consisting of twelve elements.\\
Furthermore, we directly get the following easy lemma:
\begin{lem} Fix $\alpha$ such that the moduli space of $\mathcal{U}(\alpha)$ contains a stable
representation such that $U_{i1}\subset U_{j2}$. Then the moduli spaces of $\mathcal{U}(\tau\cdot\alpha)$ and $U(\sigma_{ij}\cdot\alpha)$ contain a stable point such that $U_{j1}\subset U_{i2}$, the moduli space of $U(\sigma_{ik}\cdot\alpha)$ a stable point such that $U_{k1}\subset U_{j2}$ and finally the moduli space of $U(\sigma_{kj}\cdot\alpha)$ contains a stable representations such that $U_{i1}\subset U_{k2}$.
\end{lem}
Thus if we consider the above mentioned three cases, it suffices to restrict
to the special cases $U_{11}\subset U_{22},U_{32}$ and $U_{11}\subset U_{22},U_{21}\subset U_{12}$. We consider the first case.\\
The solutions of the systems of inequalities are determined by the methods
of Section $\ref{LU}$. In the following we denote by $P_i$ the sets of extreme points
and by $E_i$ the set of extremal rays with $i = 1, 2$. Then we have:
\[E_1=\{(1,1,1,1,1,1),(1,3,2,1,1,1),(1,2,3,1,1,1)\}\]
and
\begin{eqnarray*}
S_1&=&\{(1,1,1,0,0,0),(0,0,0,1,1,1),(0,1,0,0,0,1),(0,0,1,0,1,0),\\
&&(0,1,0,1,0,0),(0,0,1,1,0,0)\}.\end{eqnarray*}
In the following denote these vectors by $u_1, u_2, u_3$ and $w_1,\ldots,w_6$ respectively.\\
Note that, in order to calculate the extreme points and extremal rays we do
not require proper inequality. In particular, the second and third extreme
point correspond to semistable points.\\
Thus every solution $v$ of the system of inequalities is of the form
\begin{eqnarray*}
v&=&(s_1+s_2+s_3+k_1,s_1+3s_2+2s_3+k_1+k_3+k_5,\\
&&s_1+2s_2+3s_3+k_1+k_4+k_6,s_1+s_2+s_3+k_2+k_5+k_6,\\
&&s_1+s_2+s_3+k_2+k_4,s_1+s_2+s_3+k_2+k_3)\\
&=&(1+k_1,1+2s_2+s_3+k_1+k_3+k_5,1+s_2+2s_3+k_1+k_4+k_6,\\
&&1+k_2+k_5+k_6,1+k_2+k_4,1+k_2+k_3),
\end{eqnarray*}
such that $0\leq s_i \leq 1$ for $i = 1, 2, 3$ and $k_j\geq 0$ for $i\in\{1,2,3,4,5,6\}$ because
$s_1 + s_2 + s_3 = 1$. Now we are only interested in the integer-valued and also
stable solutions. We have the following:
\begin{lem}
All stable solutions $v$ with $k_i < 1$ are $v_1 =(1, 1, 1, 1, 1, 1)$,
$v_2 = (1, 2, 2, 1, 1, 1)$ and $v_3 = (1, 2, 2, 2, 1, 1)$. Otherwise we have $v = v_k +\sum_{i=1}^6 n_iw_i$ for a $n_i \in\mathbb{N}$ and a $k\in\{1, 2, 3\}$.
\end{lem}
Note that $s_i<1$ for $i\neq 1$ has to be fulfilled in order to satisfy the stability
condition because for the quivers corresponding to the extreme points except
$v_1$ there only exist semistable representations, i.e. no stable ones.\\
Furthermore, we have
\[(1,2,2,1,1,1)=\frac{1}{3}((1,1,1,1,1,1)+(1,2,3,1,1,1)+(1,3,2,1,1,1))\] and
\begin{eqnarray*}(1,2,2,2,1,1)&=&\frac{2}{3}(1,1,1,1,1,1)+\frac{1}{6}((1,2,3,1,1,1)+(1,3,2,1,1,1))\\
&&+\frac{1}{2}((0,0,1,1,0,0)+(0,1,0,1,0,0)).\end{eqnarray*}
The considered solutions have the upper bound $(1, 5, 5, 3, 2, 2)$. Thus we just
have to verify that all other solutions are linear combinations of the desired
type.\\\\
Considering the second case we obtain
\begin{eqnarray*}E_2&=&\{(1,1,1,1,1,1),(1,1,\frac{5}{2},1,1,1),(2,1,3,1,1,1),(1,2,3,1,1,1),\\&&(1,1,1,1,2,3),(1,1,1,2,1,3),(1,1,2,1,1,2),(1,1,1,1,1,\frac{5}{2})\}
\end{eqnarray*}
and
\begin{eqnarray*}
S_2&=&\{(1,1,1,0,0,0),(0,0,0,1,1,1),(1,0,0,0,0,1),(0,0,1,1,0,0),\\
&&(0,1,0,0,0,1),(0,0,1,0,1,0)\}.
\end{eqnarray*}
As above we get:
\begin{lem} All stable solutions with $k_i < 1$ are $v_1 = (1, 1, 1, 1, 1, 1)$, $v_2 =
(1, 1, 2, 1, 1, 1),$ $v_3 = (1, 1, 1, 1, 1, 2)$, otherwise we have $v = v_k +\sum_{i=1}^6n_iw_i$
for $n_i \in\mathbb{N}$, $k\in\{1,2,3\}$.
\end{lem}
Obviously we have $(1,1,2,1,1,1)=\frac{1}{3}(1,1,1,1,1,1)+\frac{2}{3}(1,1,\frac{5}{2},1,1,1)$. Again we get an upper bound, in this case $(3, 3, 5, 3, 3, 5)$. The other solutions are
again given as linear combinations.\\\\
Finally, we consider the case $\alpha_{ij} = 0$ with $i\in\{1,2,3\}$ and $j\in\{1,2\}$.
As above we can without lose of generality assume that $\alpha_{11} = 0$ that corresponds to the case $U_{11}\subset U_{22}, U_{32}$. All other inclusions can be excluded. Therefore, the extreme points and extremal rays resp. of the inequalities are given by
\begin{eqnarray*}
E_3&=&\{(0,1,1,2,2,1),(0,1,1,3,1,1),(0,1,1,1,2,2),(0,2,1,1,1,1),\\
&&(0,1,1,1,\frac{3}{2},1),(0,1,1,2,1,2),(0,1,2,1,1,1),(0,1,1,1,1,\frac{3}{2}),\\
&&(0,1,1,1,1,1)\}
\end{eqnarray*}
and
\begin{eqnarray*}
S_3&=&\{(0,0,0,1,1,1),(0,1,0,1,0,0),(0,0,1,1,0,0),(0,0,1,0,1,0),\\
&&(0,1,0,0,0,1)\}.
\end{eqnarray*}
Analogously to the other cases we obtain:
\begin{lem} All stable vectors with $k_i < 1$ are $v_1 = (0, 1, 1, 1, 1, 1)$, $v_2 =
(0, 1, 1, 2, 1, 1)$, $v_3 = (0, 2, 2, 3, 2, 2)$ or we have $v=v_k+\sum_{i=1}^6n_iw_i$ with $n_i\in\mathbb{N}$, $k\in\{1,2,3\}$.
\end{lem}
It remains to prove that we get all positive integer-valued solutions as a
unique linear combination in this way. The following lemma deals with this:
\begin{lem}
In each of the three cases the following holds: every positive
integer-valued solution is a unique linear combination of the form $v=v_k+\sum_{i=1}^6k_iw_i$ with $k_i\in\mathbb{N}$.
\end{lem}
{\it Proof.} Depending on the starting vector $v_1$ we have in each of the three
cases for a linear combination $\alpha$ that
\[\alpha_{11}+\alpha_{21}+\alpha_{31}-\alpha_{12}-\alpha_{22}-\alpha_{32}\equiv 0,1,2\modu 3.\]
This suffices to prove the uniqueness because in addition the extremal rays
are linear independent in each case. Moreover, the cases are mutually exclusive.
\qed
The next aim is to calculate the discriminant in each of these cases so that
we get a quadratic equation, whose number of integer-valued solutions determine the Euler characteristic of the considered moduli spaces.\\
Again we treat the three case from above. In the first case the solutions are
of the form
\[\alpha=(k_1+1,k_1+k_3+k_5+l,k_1+k_4+k_6+m,k_2+k_5+k_6+n,k_2+k_4+1,k_2+k_3+1)\]
with $k_i\in\mathbb{N}$ and $l = m = n = 2$ or $l = m = 2$ and $n = 1$.\\
In the second case we have the solutions
\[\alpha=(k_1+k_3+1,k_1+k_5+1,k_1+k_4+k_6+m,k_2+k_4+1,k_2+k_6+1,k_2+k_3+k_5+n)\]
with $k_i\in\mathbb{N}$ and $n = m = 1, n = 1$ and $m = 2$ or $n = 2$ and $m = 1$.\\
Finally, the solutions in the third case are given by
\[\alpha=(0,k_2+k_5+n,k_3+k_4+n,k_1+k_2+k_3+m,k_1+k_4+n,k_1+k_5+n)\]
with $k_i\in\mathbb{N}$ and $n = m = 1$ or $n = 1, 2$ and $m = n + 1$.\\\\
First we assume that $U_{i1}\nsubseteq U_{j2}$. Afterwards the discriminant in the case
$U_{i1}\subset U_{j2}$ for $i\neq j$ is obtained from this case.\\
Let $k=(k_1,k_2,k_3,k_4,k_5,k_6)$. Then we obtain in the first case
\begin{eqnarray*}
D^1(k,l,m,n)&=&2(k_1+1)^2+2(k_1+k_3+k_5+l)^2+2(k_1+k_4+k_6+m)^2\\&&+2(k_2+k_5+k_6+n)^2+2(k_2+k_4+1)^2+2(k_2+k_3+1)^2\\
&&+2(k_1+1)(k_2+k_5+k_6+n)+2(k_2+k_4+1)\\&&(k_1+k_3+k_5+l)+2(k_1+k_4+k_6+m)(k_2+k_3+1)\\&&-2(k_1+1)(k_1+k_3+k_5+l)-2(k_1+1)\\&&(k_1+k_4+k_6+m)-2(k_1+k_4+k_6+m)\\&&(k_1+k_3+k_5+l)
-2(k_2+k_5+k_6+n)(k_2+k_4+1)\\&&-2(k_2+k_5+k_6+n)(k_2+k_3+1)-2(k_2+k_4+1)\\&&(k_2+k_3+1)-4(k_1+1)(k_2+k_4+1)-4(k_1+1)\\&&(k_2+k_3+1)-4(k_1+k_3+k_5+l)(k_2+k_5+k_6+n)\\&&
-4(k_1+k_3+k_5+l)\\&&(k_2+k_3+1)-4(k_1+k_4+k_6+m)(k_2+k_4+1)\\&&-4(k_1+k_4+k_6+m)(k_2+k_3+1)
\end{eqnarray*}
\begin{eqnarray*}
\hspace{2.2cm}&=&2(k_1+1)(-k_1-3k_2-3k_3-3k_4-l-m+n-3)\\
&&+2(k_1+k_3+k_5+l)\\&&(-3k_2-k_3-k_5-3k_6+l-m-2n-1)\\&&
+2(k_1+k_4+k_6+m)\\&&(k_1-3k_2+k_3-k_4-2k_5+m-2n-1)\\&&
+2(k_2+k_5+k_6+n)(-k_2-k_3-k_4+k_5+k_6+n-2)\\&&
+2(k_2+k_4+1)(-k_3+k_4)+2(k_2+k_3+1)(k_2+k_3+1)
\end{eqnarray*}
\begin{eqnarray*}
\hspace{2.2cm}&=&-18k_1k_2-6k_1k_3-6k_1k_4-6k_1k_5-6k_1k_6-6k_2k_3\\&&-6k_2k_4-6k_2k_5
-6k_2k_6-6k_3k_5-6k_3k_6-6k_4k_5\\&&-6k_4k_6-6k_5k_6+2k_1(-3n-6)
+2k_2(-3l-3m-3)\\&&+2k_3(-3n-3)+2k_4(-3n-3)+2k_5(-3m-3)\\&&+2k_6(-3l-3)\\&&
+2(l^2+m^2+n^2-lm-2ln-2mn-2l-2m-n-2).
\end{eqnarray*}
In the second one we get
\begin{eqnarray*}
D^2(k,m,n)&=&2(k_1+k_3+1)^2+2(k_1+k_5)^2+2(k_1+k_4+k_6+m)^2\\&&+2(k_2+k_4+1)^2
2(k_2+k_6+1)^2+2(k_2+k_3+k_5+n)^2\\&&+2(k_1+k_3+1)(k_2+k_4+1)
+2(k_1+k_5+1)(k_2+k_6+1)\\&&
+2(k_1+k_4+k_6+m)(k_2+k_3+k_5+n)-2(k_1+k_3+1)\\&&(k_1+k_5+1)-2(k_1+k_3+1)(k_1+k_4+k_6+m)\\&&
-4(k_1+k_3+1)(k_2+k_6+1)-4(k_1+k_3+1)\\&&(k_2+k_3+k_5+n)
-2(k_1+k_5+1)\\&&(k_1+k_4+k_6+m)-4(k_1+k_5+1)(k_2+k_4+1)\\&&
-4(k_1+k_5+1)(k_2+k_3+k_5+n)-4(k_1+k_4+k_6+m)\\&&(k_2+k_4+1)
-4(k_1+k_4+k_6+m)(k_2+k_6+1)\\&&-2(k_2+k_4+1)(k_2+k_6+1)
2(k_2+k_4+1)\\&&(k_2+k_3+k_5+n)-2(k_2+k_6+1)(k_2+k_3+k_5+n)
\end{eqnarray*}
\vspace{-1cm}
\begin{eqnarray*}
\hspace{1.7cm}&=& 2(k_1+k_3+1)(-3k_2-k_3-3k_5-3k_6-1-m-2n)\\&&
+2(k_1+k_5+1)(-k_1-3k_2-2k_3-2k_4-k_5-m-2n)\\&&
+2(k_1+k_4+k_6+m)(k_1-3k_2-k_3-k_4-k_5-n-2)\\&&
2(k_2+k_4+1)(-k_2-k_3+k_4-k_5-k_6-n)\\&&
2(k_2+k_6+1)(-k_3-k_5-n+1)+2(k_2+k_3+k_5+n)^2
\end{eqnarray*}
\vspace{-1cm}
\begin{eqnarray*}
\hspace{1.7cm}&=&-18k_1k_2-6k_1k_3-6k_1k_4-6k_1k_5-6k_1k_6\\&&
-6k_2k_3-6k_2k_4-6k_2k_5-6k_2k_6-6k_3k_5-6k_3k_6\\&&
-6k_4k_5-6k_4k_6-12k_1-12k_2-12k_3-12k_4-12k_5\\&&-12k_6
-6k_1n-6k_2m-12m-12n+2mn+2n^2+2m^2.
\end{eqnarray*}
And finally we get
\begin{eqnarray*}
D^3(k,m,n)&=&2(k_2+k_5+n)^2+2(k_3+k_4+n)^2+2(k_1+k_2+k_3+m)^2\\&&
+2(k_1+k_4+n)^2+2(k_1+k_5+n)^2+2(k_2+k_5+n)\\&&(k_1+k_4+n)
+2(k_3+k_4+n)(k_1+k_5+n)\\&&-2(k_2+k_5+n)(k_3+k_4+n)
-2(k_1+k_2+k_3+m)\\&&(k_1+k_4+n)-2(k_1+k_2+k_3+m)(k_1+k_5+n)\\&&
-2(k_1+k_4+n)(k_1+k_5+n)-4(k_2+k_5+n)\\&&(k_1+k_2+k_3+m)
-4(k_2+k_5+n)(k_1+k_5+n)\\&&-4(k_3+k_4+n)(k_1+k_2+k_3+m)\\&&-4(k_3+k_4+n)(k_1+k_4+n)
\end{eqnarray*}
\vspace{-1cm}
\begin{eqnarray*}
\hspace{1.7cm}&=&2(k_2+k_5+n)(-3k_1-k_2-3k_3-k_5-2m-n)\\&&+2(k_3+k_4+n)(-3k_1-2k_2-k_3-k_4+k_5-2m)\\&&
+2(k_1+k_2+k_3+m)(-k_1+k_2+k_3-k_4-k_5+m-2n)\\&&
+2(k_1+k_4+n)(k_4-k_5)+2(k_1+k_5+n)(k_1+k_5+n)\end{eqnarray*}
\vspace{-1cm}
\begin{eqnarray*}
\hspace{0.3cm}&=&-6k_1k_2-6k_1k_3-6k_1k_4-6k_1k_5-6k_2k_3\\&&
-6k_2k_4-6k_2k_5-6k_3k_4-6k_3k_5-12k_1n-12k_2n\\&&-12k_3n-6k_4m-6k_5m-12mn+2m^2.
\end{eqnarray*}
Every solution of these quadratic equations with a fixed determinant $D$ automatically satisfies the inequalities. Therefore, we ''just'' have to determine
all solutions of these equations in order to calculate the Euler characteristic
of the moduli spaces of bundles of rank three. But as mentioned above we
first have to investigate the change of the discriminant in the cases in which
the subspaces are embedded in each other.\\
Obviously in all of the cases the residue class of the discriminant only depends on $l,m$ and $n$. In the first two cases we have $D\equiv 0\modu 6$ if the
starting vector is $(1, 1, 1, 1, 1, 1)$ and $D\equiv 4\modu 6$ otherwise. In the last case
we have $D\equiv 4\modu 6$ if the starting vector is $(0, 1, 1, 2, 1, 1)$ or $(0, 1, 1, 1, 1, 1)$ and $D\equiv 0\modu 6$ in the remaining case.\\
Again in the first case we have
\begin{eqnarray*}D^1_{11,22}(l,m,n,k)&=&D^1(k,l,m,n)+6\alpha_{11}\alpha_{22}\\
&=&D^1(k,l,m,n)+6(k_1k_2+k_1k_4+k_1+k_2+k_4+1)\\
&=&-6\sum_{1\leq i<j\leq 6}k_ik_j+6k_1k_4+6k_3k_4-6k_1k_2\\&&-6n(k_1+k_3+k_4)
-6l(k_2+k_6)-6m(k_2+k_5)\\&&-6(k_1+k_3+k_5+k_6)
+2(l^2+m^2\\&&+n^2-lm-2ln-2mn-2l-2m-n+1)
\end{eqnarray*}
and
\begin{eqnarray*}D^1_{11,32}(l,m,n,k)&=&D^1(k,l,m,n)+6\alpha_{11}\alpha_{32}\\
&=&D^1(k,l,m,n)+6(k_1k_2+k_1k_3+k_1+k_2+k_3+1)\\
&=&-6\sum_{1\leq i<j\leq 6}k_ik_j+6k_1k_3+6k_3k_4-6k_1k_2\\&&-6n(k_1+k_3+k_4)
-6l(k_2+k_6)-6m(k_2+k_5)\\&&-6(k_1+k_4+k_5+k_6)\\&&
+2(l^2+m^2+n^2-lm-2ln-2mn-2l\\&&-2m-n+1).
\end{eqnarray*}
In the second case we have
\begin{eqnarray*}
D^2_{11,22}(m,n,k)&=&D^2(k,m,n)+6\alpha_{11}\alpha_{22}\\
&=&-6\sum_{1\leq i<j\leq 6}k_ik_j-6k_1k_2+6k_1k_6+6k_2k_3+6k_3k_4\\&&+6k_3k_6+6k_5k_6
-6(k_1+k_2+k_3+2k_4+2k_5+k_6)\\&&-6k_1n-6k_2m-12m-12n+2mn+2n^2+2m^2+6
\end{eqnarray*}
and
\begin{eqnarray*}
D^2_{21,12}(m,n,k)&=&D^2(k,m,n)+6\alpha_{21}\alpha_{12}\\
&=&-6\sum_{1\leq i<j\leq 6}k_ik_j-6k_1k_2+6k_1k_4+6k_2k_5+6k_3k_4\\&&+6k_4k_5+6k_5k_6
-6(k_1+k_2+2k_3+k_4+k_5+2k_6)\\&&-6k_1n-6k_2m-12m-12n+2mn+2n^2+2m^2+6.
\end{eqnarray*}
In the third case the discriminant stays constant because $\alpha_{11} = 0$.\\
We only have to consider the last four quadratic equations and the solutions
in the third case, i.e. $\alpha_{11} = 0$, because the moduli spaces are $\mathbb{P}^1$ without
two points. In particular, the Euler characteristic is zero.\\
If we evaluate the above functions at the relevant points for $m, n$ and $l$, we
obtain:
\begin{eqnarray*}D^1_{11,22}(2,2,1,k)&=&-6\sum_{1\leq i<j\leq 6}k_ik_j+6k_1k_4+6k_3k_4-6k_1k_2\\&&-6(2k_1+4k_2+2k_
3+k_4+3k_5+3k_6)-22
\end{eqnarray*}
and
\begin{eqnarray*}D^1_{11,22}(2,2,2,k)&=&-6\sum_{1\leq i<j\leq 6}k_ik_j+6k_1k_4+6k_3k_4-6k_1k_2\\&&-6(3k_1+4k_2+3k_
3+2k_4+3k_5+3k_6)-34
\end{eqnarray*}
respectively, and
\begin{eqnarray*}D^1_{11,32}(2,2,1,k)&=&-6\sum_{1\leq i<j\leq 6}k_ik_j+6k_1k_3+6k_3k_4-6k_1k_2\\&&-6(2k_1+4k_2+k_
3+2k_4+3k_5+3k_6)-22
\end{eqnarray*}
and
\begin{eqnarray*}D^1_{11,32}(2,2,2,k)&=&-6\sum_{1\leq i<j\leq 6}k_ik_j+6k_1k_3+6k_3k_4-6k_1k_2\\&&-6(3k_1+4k_2+2k_
3+3k_4+3k_5+3k_6)-34.
\end{eqnarray*}
respectively. Evaluating $D^2$, we get
\begin{eqnarray*}
D^2_{11,22}(1,2,k)&=&-6\sum_{1\leq i<j\leq 6}k_ik_j-6k_1k_2+6k_1k_6+6k_2k_3+6k_3k_4\\&&+6k_3k_6+6k_5k_6
-6(3k_1+2k_2+k_3+2k_4+2k_5+k_6)\\&&-16
\end{eqnarray*}
and
\begin{eqnarray*}
D^2_{11,22}(2,1,k)&=&-6\sum_{1\leq i<j\leq 6}k_ik_j-6k_1k_2+6k_1k_6+6k_2k_3+6k_3k_4\\&&+6k_3k_6+6k_5k_6
-6(2k_1+3k_2+k_3+2k_4+2k_5+k_6)\\&&-16
\end{eqnarray*}
respectively, and
\begin{eqnarray*}
D^2_{21,12}(1,2,k)&=&-6\sum_{1\leq i<j\leq 6}k_ik_j-6k_1k_2+6k_1k_4+6k_2k_5+6k_3k_4\\&&+6k_4k_5+6k_5k_6
-6(3k_1+2k_2+2k_3+k_4+k_5+2k_6)\\&&-16
\end{eqnarray*}
and
\begin{eqnarray*}
D^2_{21,12}(2,1,k)&=&-6\sum_{1\leq i<j\leq 6}k_ik_j-6k_1k_2+6k_1k_4+6k_2k_5+6k_3k_4\\&&+6k_4k_5+6k_5k_6
-6(2k_1+3k_2+2k_3+k_4+k_5+2k_6)\\&&-16
\end{eqnarray*}
respectively. Further if $k=(k_1,k_2,k_3,k_4,k_5)$ we get
\[D^3(2,1,k)=-6\sum_{1\leq i<j\leq 5}k_ik_j+6k_4k_5
-12(k_1+k_2+k_3+k_4+k_5)-16\]
and
\[D^3(1,1,k)=-6\sum_{1\leq i<j\leq 5}k_ik_j+6k_4k_5
-6(2k_1+2k_2+2k_3+k_4+k_5)-10\]
respectively. Obviously we have the following equations:
\[D^1_{11,22}(k_1,k_2,k_3,k_4,k_5,k_6)=D^1_{11,32}(k_1,k_2,k_4,k_3,k_5,k_6),\]
\[D^2_{11,22}(k_1,k_2,k_3,k_4,k_5,k_6)=D^2_{21,12}(k_1,k_2,k_5,k_6,k_3,k_4)\]
and
\[D^2_{11,22}(1,2,k_1,k_2,k_3,k_4,k_5,k_6)=D^2_{11,22}(2,1,k_2,k_1,k_6,k_5,k_4,k_3).\]
Define
\[K^1_{i,j}(D,l,m,n)=\{k=(k_1,k_2,k_3,k_4,k_5,k_6)\in\mathbb{N}^6_0\mid D^1_{i,j}(l,m,n,k)=D\}\]
for suitable $i,j$. Define $K^2_{i,j}(D,m,n)$ and $K^3(D, m, n)$ analogously.\\
By use of $\ref{six}$, the preceding calculations and the mentioned equalities we get:
\begin{satz} Let $D\equiv 4\modu 6$. Then we have
\begin{eqnarray*} \chi(\mathcal{M}(D))&=&6|K^1_{11,22}(D,2,2,1)|+6|K^1_{11,22}(D,2,2,2)|\\&&+6|K^2_{11,22}(D,1,2)|+3|K^3(D,2,1)|+3|K^3(D,1,1)|.\end{eqnarray*}
\end{satz}
Analogously to the case of rank two bundles, we obtain the following corollary concerning the generating function of the Euler characteristic:
\begin{kor} Let $D\equiv 4\modu 6$. We have
\begin{eqnarray*}F(x)&=&6(\sum_{k\in\mathbb{N}^6_0}x^{D^1_{11,22}(2,2,1,k)}+x^{D^1_{11,22}(2,2,2,k)}+x^{D^2_{11,22}(1,2,k)})\\&&+3\sum_{k\in\mathbb{N}_0^5} x^{D^3(2,1,k)}+3\sum_{k\in\mathbb{N}_0^5} x^{D^3(1,1,k)}.
\end{eqnarray*}
\end{kor}
Thus we have \begin{eqnarray*}
F(x)&=&0x^{-4}+3x^{-10}+15x^{-16}+36x^{-22}+69x^{-28}+114x^{-34}+165x^{-40}\\&&+246x^{-46}+303x^{-52}+432x^{-58}+492x^{-64}+669x^{-70}+726x^{-76}\\&&
+975x^{-82}+999x^{-88}+1332x^{-94}+1338x^{-100}+1743x^{-106}\\&&+1716x^{-112}+2226x^{-118}+2130x^{-124}+2775x^{-130}+2625x^{-136}\\&&+3354x^{-142}
+3129x^{-148}+4041x^{-154}+3735x^{-160}+4752x^{-166}\\&&+4317x^{-172}+5532x^{-178}+5070x^{-184}+6393x^{-190}+\mathcal{O}(x^{-202}).
\end{eqnarray*}
\subsection{The case $D\equiv 0\modu 6$}
In this section we discuss the case of the discriminants satisfying $D\equiv 0\modu 6$. The main difference to the preceding case is that there also exist semistable points. Thus we have to modify the methods slightly.\\
First let $\alpha_{ij} \neq 0$ and consider all inclusions of vector spaces pointed out
in the last section. All extreme points except the point $(1, 1, 1, 1, 1, 1)$ correspond to points in the case $D\equiv 4\modu 6$. Thus it remains to consider this extreme point. But we have to keep in mind that the inequalities can be satisfied with equality. Thus the inequalities do not exclude each other.\\
Therefore, we consider all extremal rays appearing in the last section, i.e.
\begin{eqnarray*}S&=&\{(1,1,1,0,0,0),(0,0,0,1,1,1),(1,0,0,0,1,0),(0,1,0,0,0,1),\\  
&&(0,0,1,1,0,0),(1,0,0,0,0,1),(0,1,0,1,0,0),(0,0,1,0,1,0)\}.
\end{eqnarray*}
Note that the extremal rays in the nine cases of the last section arise from
these eight rays by considering the first two and in addition removing one
of the rays three to five and one of the rays six to eight. The linear combinations having as starting point the extreme point $(1, 1, 1, 1, 1, 1)$ thus have
the following standard form
\begin{eqnarray*}\alpha&=&(k_1+k_3+k_6+1,k_1+k_4+k_7+1,\\&&k_1+k_5+k_8+1,k_2+k_5+k_7+1,k_2+k_3+k_8+1,\\&&k_2+k_4+k_6+1).\end{eqnarray*}
Let $k=(k_1,k_2,\ldots,k_8)$. Then the discriminant is given by
\begin{eqnarray*}
D(k)&=&2(k_1+k_3+k_6+1)(k_1+k_3+k_6+1)+2(k_1+k_3+k_6+1)\\&&(k_2+k_5+k_7+1)-2(k_1+k_3+k_6+1)(k_1+k_4+k_7+1)\\&&
-2(k_1+k_3+k_6+1)(k_1+k_5+k_8)-4(k_1+k_3+k_6+1)\\&&(k_2+k_3+k_8)-4(k_1+k_3+k_6)(k_2+k_4+k_6+1)\\&&
+2(k_1+k_4+k_7+1)(k_1+k_4+k_7+1)+2(k_1+k_4+k_7+1)\\&&(k_2+k_3+k_8+1)-2(k_1+k_4+k_7+1)(k_1+k_5+k_8+1)\\&&
-4(k_1+k_4+k_7+1)(k_2+k_5+k_7+1)-4(k_1+k_4+k_7+1)\\&&(k_2+k_4+k_6+1)+2(k_1+k_5+k_8+1)(k_1+k_5+k_8+1)\\&&
+2(k_1+k_5+k_8+1)(k_2+k_4+k_6+1)-4(k_1+k_5+k_8+1)\\&&(k_2+k_5+k_7+1)-4(k_1+k_5+k_8+1)(k_2+k_3+k_8+1)\\&&
+2(k_2+k_5+k_7+1)(k_2+k_5+k_7+1)-(k_2+k_5+k_7+1)\\&&2(k_2+k_3+k_8+1)-2(k_2+k_5+k_7+1)(k_2+k_4+k_6+1)\\&&
+2(k_2+k_3+k_8+1)(k_2+k_3+k_8+1)-2(k_2+k_3+k_8+1)\\&&(k_2+k_4+k_6+1)
+2(k_2+k_4+k_6+1)(k_2+k_4+k_6+1)
\end{eqnarray*}
\vspace{-1cm}
\begin{eqnarray*}
\hspace{1.1cm}&=&2(k_1+k_3+k_6+1)(-k_1-3k_2-k_3-3k_4-k_6-3k_8-4)\\&&
+2(k_1+k_4+k_7+1)(-3k_2+k_3-k_4-3k_5-2k_6-k_7-3)\\&&
+2(k_1+k_5+k_8+1)(k_1-3k_2-2k_3+k_4-k_5+k_6-2k_7\\&&-k_8-2)
+2(k_2+k_5+k_7+1)(-k_2-k_3-k_4+k_5-k_6\\&&+k_7-k_8-1)
+2(k_2+k_3+k_8+1)(k_3-k_4-k_6+k_8)\\&&+2(k_2+k_4+k_6+1)(k_2+k_4+k_6+1)
\end{eqnarray*}
\vspace{-1cm}
\begin{eqnarray*}
\hspace{1.1cm}&=&-6\sum_{1\leq i<j\leq 8}k_ik_j-12k_1k_2+6k_3k_7+6k_4k_8+6k_5k_6-18k_1\\&&-18k_2-12\sum_{i=3}^8k_i-18
\end{eqnarray*}
Define $K_1 = \{k_3, k_4, k_5\}$ and $K_2 = \{k_6, k_7, k_8\}$. Again consider a filtration
in standard form:
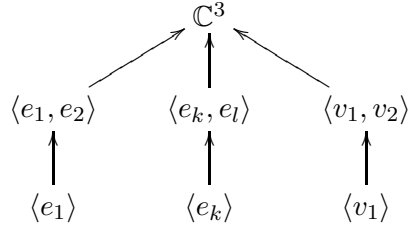
\filt{e_1}{e_2}{e_k}{e_l}{v_1}{v_2}
Thereby $k, l\in\{1,2,3\}$ such that $k\neq l$ and $v_1\neq v_2$ are arbitrary vectors.\\
The stable filtrations are filtrations of the form
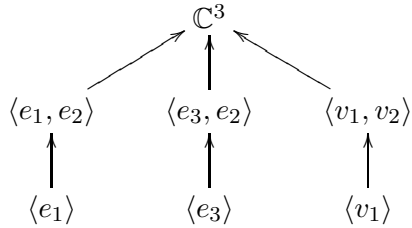
\filt{e_1}{e_2}{e_3}{e_2}{v_1}{v_2}
with certain conditions for $v_1$ and $v_2$ investigated in more detail now.\\
Define $e_{ij} = e_i + e_j$ and $e_{123} = (1, 1, 1)$. If $v_1 = e_i$, all resulting representations are unstable. If $v_1 = e_{ij}$ for $i\neq j$, the obtained filtrations correspond to polystable representations that will be analysed later. Thus let $(v_1)_i\neq 0$
and we may without lose of generality assume that $v_1 = e_{123}$.\\
Furthermore, we may assume that $(v_2)_1 = 0$. Indeed, we can add arbitrary
multiples of $e_{123}$ to $v_2$.
\begin{lem}\label{koecher} Every stable representation of the quiver $\mathcal{U}(1,1,1,1,1,1)$  is
given by
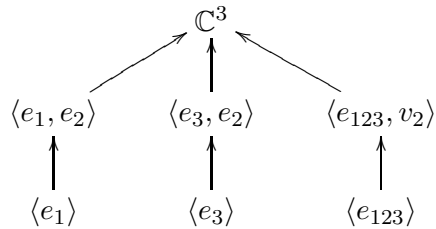
\filt{e_1}{e_2}{e_3}{e_2}{e_{123}}{v_2}
where $(v_2)_1 = 0$ and $v_2 \neq e_2, e_3, e_{23}$. In particular, we have for the moduli
space of stable representations $\mathcal{M}(1,1,1,1,1,1)^s=\mathbb{P}^1\backslash\{(1:0),(0:1),(1:1)\}$.
\end{lem}
Note that the case $v_2 = e_{23}$ is equivalent to the case $v_2 = e_1$. It is seen by an
easy calculation that the three filtrations given by the points $(1 : 0), (0 : 1)$
and $(1 : 1)$ are not stable but semistable.\\
We obtain the following corollary:
\begin{kor} 
If there exist stable points for $\mathcal{U}(\alpha_{11},\alpha_{12},\alpha_{21},\alpha_{22},\alpha_{31},\alpha_{32})$,
the points corresponding to the points of Lemma $\ref{koecher}$ are already stable.
\end{kor}
The polystable points can be described as follows: for $\alpha_{ij} = k$ for all $i, j$ we
obtain a polystable representation induced by the following:
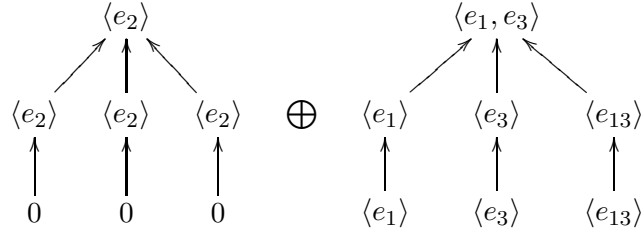
\[
\begin{xy}
\xymatrix@R20pt@C10pt{
&\langle e_2\rangle &&&&\langle e_1,e_3\rangle&\\\langle e_2\rangle\ar[ru]&\langle e_2\rangle\ar[u]&\langle e_2\rangle\ar[lu]&\bigoplus&\langle e_1\rangle\ar[ru]&\langle e_3\rangle\ar[u]&\langle e_{13}\rangle\ar[lu]\\0\ar[u]&0\ar[u]&0\ar[u]&&\langle e_1\rangle\ar[u]&\langle e_3\rangle\ar[u]&\langle e_{13}\rangle\ar[u]}
\end{xy}
\]
This point also induces a polystable point of the quiver obtained by extending the arms with the vector $(1, 1, 1, 1, 1, 1)$. Call these points polystable of
type $1$.\\
Further consider the polystable representations
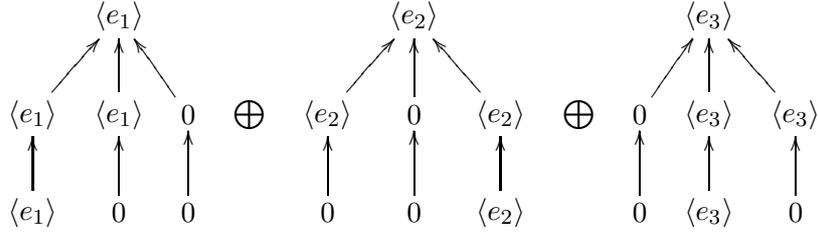
\[
\begin{xy}
\xymatrix@R20pt@C8pt{
&\langle e_1\rangle &&&&\langle e_2\rangle &&&&\langle e_3\rangle &\\\langle e_1\rangle\ar[ru]&\langle e_1\rangle\ar[u]&0\ar[lu]&\bigoplus&\langle e_2\rangle\ar[ru]&0\ar[u]&\langle e_2\rangle\ar[lu]&\bigoplus&0\ar[ru]&\langle e_3\rangle\ar[u]&\langle e_3\rangle\ar[lu]\\\langle e_1\rangle\ar[u]&0\ar[u]&0\ar[u]&&0\ar[u]&0\ar[u]&\langle e_2\rangle\ar[u]&&0\ar[u]&\langle e_3\rangle\ar[u]&0\ar[u]}
\end{xy}
\]
and
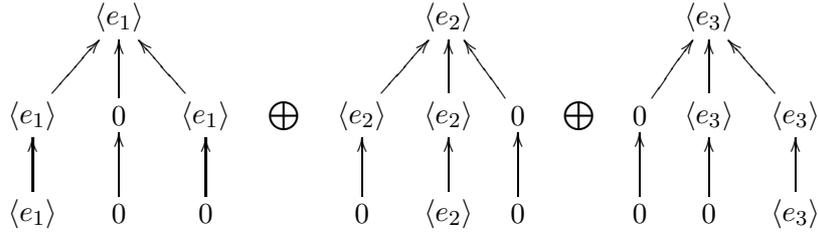
\[
\begin{xy}
\xymatrix@R20pt@C8pt{
&\langle e_1\rangle &&&&\langle e_2\rangle &&&&\langle e_3\rangle &\\\langle e_1\rangle\ar[ru]&0\ar[u]&\langle e_1\rangle\ar[lu]&\bigoplus&\langle e_2\rangle\ar[ru]&\langle e_2\rangle\ar[u]&0\ar[lu]&\bigoplus&0\ar[ru]&\langle e_3\rangle\ar[u]&\langle e_3\rangle\ar[lu]\\\langle e_1\rangle\ar[u]&0\ar[u]&0\ar[u]&&0\ar[u]&\langle e_2\rangle\ar[u]&0\ar[u]&&0\ar[u]&0\ar[u]&\langle e_3\rangle\ar[u]}
\end{xy}
\]
These are the remaining polystable points in the case $\alpha_{ij} = 1$, which per\-sist under the extensions given by $(\alpha_{21}, \alpha_{12})$, $(\alpha_{31}, \alpha_{22})$ and $(\alpha_{11}, \alpha_{32})$ and $(\alpha_{11}, \alpha_{22})$, $(\alpha_{21}, \alpha_{32})$ and $(\alpha_{31}, \alpha_{12})$ respectively. Call these points polystable
of type 2.\\
If we consider the case with lengths of arms given by $(1, 1, 2, 2, 1, 1)$, we
get the polystable point
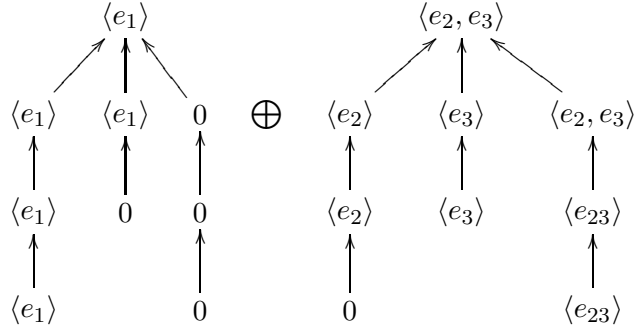
\[
\begin{xy}
\xymatrix@R20pt@C10pt{
&\langle e_1\rangle &&&&\langle e_2,e_3\rangle&\\\langle e_1\rangle\ar[ru]&\langle e_1\rangle\ar[u]&0\ar[lu]&\bigoplus&\langle e_2\rangle\ar[ru]&\langle e_3\rangle\ar[u]&\langle e_2,e_3\rangle\ar[lu]\\\langle e_1\rangle\ar[u]&0\ar[u]&0\ar[u]&&\langle e_2\rangle\ar[u]&\langle e_3\rangle\ar[u]&\langle e_{23}\rangle\ar[u]\\\langle e_1\rangle\ar[u]&&0\ar[u]&&0\ar[u]&&\langle e_{23}\rangle\ar[u]}
\end{xy}
\]
and the corresponding polystable points according to the remaining five
extensions.\\
They persist under extensions given by $(\alpha_{11}, \alpha_{32})$, $(\alpha_{31}, \alpha_{22})$, $(\alpha_{21}, \alpha_{12})$ and $(\alpha_{31}, \alpha_{12})$. Call these points polystable of type 3. Note that polystable points of type 2 and 3 also persist under extending the arms by the vectors $(1, 1, 1, 0, 0, 0)$ and $(0, 0, 0, 1, 1, 1)$ respectively.\\
If we want to determine a solution of the given system of inequalities, we
may assume that there exists at least one $k\in K_i$ for every $i = 1, 2$ with
$k = 0$. If we again consider the inequalities ($\ref{eins}$)-($\ref{vier}$) and the investigations of the polystable points, we get in conclusion:
\begin{lem}
Let $\alpha\in\mathbb{N}^6_0$ be in standard form.
\begin{enumerate}
\item The filtration induced by \filt{e_1}{e_2}{e_3}{e_2}{e_{123}}{e_2}
is stable if and only if $k_1 > k_2$.
\item The filtration induced by
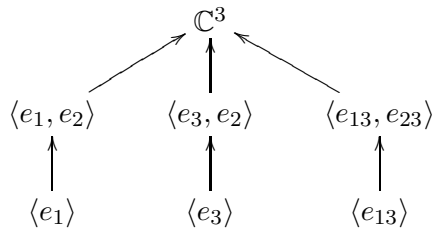
\filt{e_1}{e_2}{e_3}{e_2}{e_{13}}{e_{23}}
is stable if and only if $k_2 > k_1$.
\item If $k_1 = k_2$, there exists exactly one polystable point of type 1.
\item If $k\neq 0$ for exactly one $k\in K_i$, $i = 1, 2$, there exists a polystable point
of type 2 and a polystable point of type 3.
\item If $k, l\neq 0$ for exactly two different $k, l\in K_i$, $i = 1, 2$, there exists a
stable point such that $U_{i1} \subset U_{j2}$ with $i\neq j$ and a polystable point of
type 2.
\item If $k, l\neq 0$ for exactly one $k\in K_1$ and exactly one $l\in K_2$, there exist
two polystable points of type 3.
\item If $k, l, n\neq 0$ for exactly two different $k, l\in K_i$ and one $n\in K_j$ with
$i\neq j$, there exists a stable point such that $U_{i1}\subset U_{j2}$ with $i\neq j$ and a
polystable point of type 3.
\item If $k, l, m, n\neq 0$ for exactly two different $k, l\in K_1$ and exactly two
different $n,m\in K_2$, there exist two stable points such that $U_{i1} \subset U_{j2}$ with $i\neq j$.
\end{enumerate}
\end{lem}
Considering all polystable filtrations treated in this section, an easy calculation using the results of $\cite{alb}$ shows that all of these points are smooth. If
we assume that there exists a least one stable filtration, in the same manner
as in the last section we get that
\[M(\mathcal{U}(\alpha))^{ss}\cong\mathbb{P}^1.\]
The moduli space of stable points is obtained by the considerations of the
last lemma.\\
Note that the moduli space of stable points and the one of semistable points
coincide if in the eighth case of the lemma $k_1\neq k_2$ holds.\\
Finally, we have to consider the cases when $\alpha_{ij} = 0$ for exactly one pair $i, j$.
This is the case of the extreme point $(0, 2, 2, 3, 2, 2)$. As in the last section,
we may assume $\alpha_{11} = 0$. Then the discriminant is given by
\begin{eqnarray*}D^0(k)&=&-6\sum_{1\leq i<j\leq 5}k_ik_j+6k_4k_5-24k_1-24k_2\\
&&-24k_3-18k_4-18k_5-54.
\end{eqnarray*} 
Again we determine the generating function. Therefore, define
\[D^1(k)=D(k)+6(k_1+k_5+k_8+1)(k_2+k_5+k_7+1),\]
\[D^2(k)=D(k)+6(k_1+k_3+k_6+1)(k_2+k_4+k_6+1)\]
and
\[D^3(k)=D(k)+6(k_1+k_5+k_8+1)(k_2+k_3+k_8+1).\]
This corresponds to the discriminant in the case of inclusions $U_{31}\subset U_{12}$,
 $U_{11}\subset U_{32}$ and $U_{31}\subset U_{22}$. It suffices to consider the cases provided by this three cases, the other ones can again be constructed via permutation of the
arms.\\
Let
\[F(x)=\sum_{n=0}^{\infty}\chi(\mathcal{M}(6n))x^{-6n}\]
be the generating function of the Euler characteristic.
Define $\mathcal{N}^{k,l}=\mathbb{N}^k\times \mathbb{N}_+^{l}$ for a $k\in\mathbb{N}_+$.\\
Then we get the following result:
\begin{satz} We have
\begin{eqnarray*}
F(x)&=&-\sum_{k=0}^{\infty}x^{D(k,k,0,\ldots,0)}-6\sum_{k\in\mathcal{N}^{1,1}}x^{D(k_1,k_1,k_2,0,\ldots,0)}\\
&&-3\sum_{k\in\mathcal{N}^{1,2}}x^{D(k_1,k_1,k_2,0,0,0,k_3,0)}-6\sum_{k\in\mathcal{N}^{1,2}}x^{D(k_1,k_1,k_2,0,0,0,0,k_3)}\\
&&-6\sum_{k\in\mathcal{N}^{1,2}}x^{D(k_1,k_1,k_2,k_3,0,0,0,0)}-12\sum_{k\in\mathcal{N}^{1,3}}x^{D(k_1,k_1,k_2,k_3,0,0,k_4,0)}\\
&&-6\sum_{k\in\mathcal{N}^{1,3}}x^{D(k_1,k_1,k_2,k_3,0,k_4,0,0)}-3\sum_{k\in\mathcal{N}^{1,4}}x^{D(k_1,k_1,k_2,k_3,0,0,k_4,k_5)}\\
&&-6\sum_{k\in\mathcal{N}^{1,4}}x^{D(k_1,k_1,k_2,k_3,0,k_4,k_5,0)}\\
&&+6\sum_{k\in\mathcal{N}^{2,2}}x^{D^1(k_1,k_2,k_3,k_4,0,0,0,0)}+12\sum_{k\in\mathcal{N}^{2,3}}x^{D^1(k_1,k_2,k_3,k_4,0,0,k_5,0)}\\
&&+6\sum_{k\in\mathcal{N}^{2,3}}x^{D^1(k_1,k_2,k_3,k_4,0,k_5,0,0)}+3\sum_{k\in\mathcal{N}^{2,4}}x^{D^1(k_1,k_2,k_3,k_4,0,0,k_5,k_6)}\\
&&+3\sum_{k\in\mathcal{N}^{2,4}}x^{D^2(k_1,k_2,k_3,k_4,0,0,k_5,k_6)}+6\sum_{k\in\mathcal{N}^{2,4}}x^{D^1(k_1,k_2,k_3,k_4,0,k_5,k_6,0)}\\
&&+6\sum_{k\in\mathcal{N}^{2,4}}x^{D^3(k_1,k_2,k_3,k_4,0,k_5,k_6,0)}+6\sum_{k\in\mathbb{N}^5}x^{D^0(k)}.
\end{eqnarray*}
\end{satz}
As far as the case $k_1 \neq k_2$ is concerned, note that we only have to count
the stable points coming from inclusions. Indeed, the moduli space is the
projective line without two points.\\
Thus we have \begin{eqnarray*}
F(x)&=&0x^0+0x^{-6}+0x^{-12}-1x^{-18}+0x^{-24}-6x^{-30}+0x^{-36}-3x^{-42}\\&&
-12x^{-48}+6x^{-54}+12x^{-60}-15x^{-66}+17x^{-72}+72x^{-78}-24x^{-84}\\&&+102x^{-90}+30x^{-96}+138x^{-102}+132x^{-108}+171x^{-114}+27x^{-120}\\&&+420x^{-126}+204x^{-132}+360x^{-138}+180x^{-144}+678x^{-150}\\&&+192x^{-156}+773x^{-162}+351x^{-168}+906x^{-174}+624x^{-180}\\&&+816x^{-186}+519x^{-192}+\mathcal{O}(x^{-198}).
\end{eqnarray*}

\end{document}